\documentclass[a4paper,twoside,12pt]{article}

\usepackage{amsmath}
\usepackage{amsthm}
\usepackage{amsfonts}
\usepackage{amssymb, amscd}

 \numberwithin{equation}{section}
\title{Deformations of Lie algebras using $\s$-derivations}
\author{Jonas T. Hartwig\\
   \footnotesize \emph{Department of Mathematics,
    Chalmers University of Technology
    and G\"oteborg}\\ \footnotesize \emph{University
    SE-412 96 G\"oteborg, Sweden,
    hart@math.chalmers.se} \normalsize
    \and
    Daniel Larsson, Sergei Silvestrov\\
   \footnotesize \emph{Centre for Mathematical Sciences, Department of
   Mathematics, Lund Institute}\\
   \footnotesize \emph{of Technology, Lund University, Box
   118, SE-221 00 Lund, Sweden}\\
   \footnotesize \emph{dlarsson@maths.lth.se,  sergei.silvestrov@math.lth.se}\normalsize\and
}

\date{November 12, 2003}
\theoremstyle{definition}
\newtheorem{dfn}{Definition}
\newtheorem{example}{Example}
\newtheorem{remark}{Remark}

\theoremstyle{plain}
\newtheorem{prop}{Proposition}
\newtheorem{thm}[prop]{Theorem}

\newtheorem{lem}[prop]{Lemma}

\def\<{ \langle }
\def\>{ \rangle }
\def\A{ \mathcal{A} }

\def\a{ \alpha }

\def\c{ \mathbf{c} }

\def\d{ \delta }

\def\s{ \sigma }
\def\vs{ \varsigma }
\def\t{ \tau }

\def\dh{\tilde d}
\def\C{ \mathbb{C} }
\def\Z{ \mathbb{Z} }
\def\N{ \mathbb{N} }

\def\De{ \Delta }
\def\hvs{ \hat{\varsigma} }
\def\mfD{ \mathfrak{D} }
\def\M{ \mathfrak{M} }
\def\Vir{ \mathrm{Vir} }

\def\cs{ \circlearrowleft }

\def\afrak{ \mathfrak{a} }

\def\kbf{ \mathbf{k} }
\def\lbf{ \mathbf{l} }
\def\hbf{ \mathbf{h} }

 \DeclareMathOperator{\id}{id}

\DeclareMathOperator{\Ann}{Ann}
\DeclareMathOperator{\GCD}{gcd}\DeclareMathOperator{\pr}{pr}\DeclareMathOperator{\sign}{sign}

\begin{document}

\maketitle

\begin{abstract}
In this article we develop an approach to deformations of the Witt
and Virasoro algebras based on $\sigma$-derivations. We show that
$\sigma$-twisted Jacobi type identity holds for generators of such
deformations. For the $\sigma$-twisted generalization of Lie
algebras modeled by this construction, we develop a theory of
central extensions. We show that our approach can be used to
construct new deformations of Lie algebras and their central
extensions, which in particular include naturally the
$q$-deformations of the Witt and Virasoro algebras associated to
$q$-difference operators, providing also corresponding
$q$-deformed Jacobi identities. 
\end{abstract}

\footnotetext[2]{The research was supported by the
Crafoord Foundation, the Swedish Royal Academy of Sciences and
Mittag-Leffler Institute. The results of this paper have been
reported at the Non-commutative Geometry workshop, Mittag-Leffler
Institute, Stockholm, September 8-12, 2003.}

\newpage
\tableofcontents
\section{Introduction}
Multiplicative deformations-discretizations of derivatives have
many applications in models of quantum phenomena, as well as in
analysis of complex systems and processes exhibiting complete or
partial scaling invariance. The key algebraic property which is
shared by these differential and difference type operators and
making them so useful is that they satisfy some versions of the
Leibniz rule explaining how to calculate the operator on products
given its action on each factor. It is desirable therefore to have
a single unifying differentiation theory, which would be concerned
with operators of a certain general class, satisfying generalized
Leibniz rule and containing as examples the classical
differentiation and other well-known derivations and differences.

The infinite-dimensional Lie algebra of complex polynomial vector
fields on the unit circle, the Witt algebra, is an important
example in the classical differential and integral calculus,
relating it to topology and geometry, and at the same time
responsible for many of its key algebraic properties. The
universal enveloping algebra of the Witt algebra is isomorphic to
an associative algebra with an infinite number of generators
$\{d_j \mid j \in \mathbb{Z} \}$ and defining relations
\begin{equation} \label{rel:Witt}
[d_n,d_m]=d_n d_m - d_n d_m = (n-m) d_{n+m} \quad \mbox{ for } n,m
\in \mathbb{Z}.
\end{equation}

The Witt algebra can also be defined  as the complex Lie algebra
of derivations on the algebra of Laurent polynomials $\mathbb{C}
[t,t^{-1}]$ in one variable, that is the Lie algebra of linear
operators $D$ on $\mathbb{C} [t,t^{-1}]$ satisfying the ordinary
Leibniz rule $D(ab) = D(a) b + a D(b)$, with commutator taken as
the Lie algebra product. This definition will be most important in
this article, as it will be taken as a starting point for
generalization of the Witt algebra, incorporating operators
obeying a generalized Leibniz rule twisted by an endomorphism
(Definition \ref{def:genWittalg}).

Important examples of such twisted derivation-type operators, extensively investigated in physics and
engineering and lying at the foundations of \mbox{$q$-analysis,} are the Jackson $q$-derivative $D_q (f)(t) =
\frac{f(qt) - f(t)}{qt-t}$ and $q$-derivative $M_t D_q (f) (t) = \frac{f(qt) - f(t)}{q-1}$ acting on $\mathbb{C}
[t,t^{-1}]$ or various function spaces. It satisfies a $\sigma_q$-twisted ($q$-deformed) Leibniz rule $D(fg) =
D(f) g + \sigma_q(f) D(g)$ for the re-scaling automorphism $\sigma_q (f) (t)= f(qt)$. In this special case our
general construction yields a natural $q$-deformation of the Witt algebra which becomes the usual Witt algebra
defined by \eqref{rel:Witt} when $q=1$ (Theorem \ref{th:qWittlinear}). It is closely related to $q$-deformations
of the Witt algebra introduced and studied in \cite{AizSato1,ChaiIsLukPopPres,ChaiKuLuk,ChaiPres3,
ChungWS,CurtrZachos1,Kassel1,Polychronakos,SatoHT1,SatoHT2}. However, our defining commutation relations in this
case look somewhat different, as we obtained them, not from some conditions aiming to resolve specifically the
case of $q$-deformations, but rather by choosing $\mathbb{C} [t,t^{-1}]$ as an example of the underlying
coefficient algebra and specifying $\sigma$ to be the automorphism $\sigma_q$ in our general construction for
$\sigma$-derivations. By simply choosing a different coefficient algebra or basic $\sigma$-derivation one can
construct many other analogues and deformations of the Witt algebra. We demonstrate this by examples,
constructing a class of deformations of the Witt algebra parameterized by integers defining arbitrary
endomorphisms of $\mathbb{C} [t,t^{-1}]$ (Theorem \ref{th:qWittnonlinear}), and also constructing a
multi-dimensional analogue of the Witt algebra deforming commutative algebra with infinite countable number of
generators, by taking the underlying algebra to be Laurent polynomials in several variables $\C[z_1^{\pm
1},z_2^{\pm 1},\dots,z_n^{\pm 1}]$ and choosing $\s$ to be the endomorphism mapping $z_1,\dots,z_n$ into
monomials (Theorem \ref{th:qWittmultdimmonom}). The important feature of our approach is that, as in the
non-deformed case, the deformations and analogues of Witt algebra obtained by various choices of the underlying
coefficient algebra, of the endomorphism $\s$ and of the basic $\s$-derivation, are precisely the natural
algebraic structures for the differential and integral type calculi and geometry based on the corresponding
classes of generalized derivation and difference type operators.

The non-deformed Witt algebra has a unique, up to multiplication
by a scalar, one-dimensional central Lie algebra extension, the
Virasoro Lie algebra. Its universal enveloping algebra, also
usually called the Virasoro algebra, is the algebra with infinite
set of generators $\{d_j \mid j \in \mathbb{Z} \} \cup
\{\mathbf{c}\}$ and defining relations
\begin{equation} \label{rel:Virasoro}
\begin{array}{ll}
& [d_j,d_k] = d_j d_k - d_k d_j = (j-k) d_{j+k} + \delta_{j+k,0}
\frac{1}{12} (j+1)j(j-1) \c, \\
& [ \c,d_k ] = \c d_k - d_k \c = 0, \quad \mbox{ for } j,k \in
\mathbb{Z}.
\end{array}
\end{equation}
We develop in this article a framework for construction of central extensions of deformed Witt algebras built on
$\sigma$-derivations. To this end we show first that our generalization of the Witt algebra to general
$\s$-derivations satisfies skew-symmetry and a generalized (twisted) Jacobi identity (Theorem
\ref{thm:GenWitt}). The generalized Jacobi identity \eqref{eq:GenWittJacobi} has six terms, three of them
twisted from inside and the other three twisted on the outside. This defines a class of non-associative algebras
with multiplication satisfying skew-symmetry and such generalized Jacobi identities, and containing Lie algebras
as the untwisted case. Sometimes the twisting can be put on the inside of all terms of the generalized Jacobi
identity in the same way, and  the terms can be coupled to yield the generalized Jacobi identity with three
terms. For example, this is the case for the $q$-deformation of the Witt algebra in Theorem
\ref{th:qWittlinear}. Armed with this observation we define the corresponding class of non-associative algebras,
calling it in this article hom-Lie algebras (Definition \ref{def:homLiealg}, Section \ref{subsec:homLiealg}),
since it is associated with a twisting homomorphism. When the twisting homomorphism is the identity map, the
generalized Jacobi identity becomes twice the usual Jacobi identity for Lie algebras, making Lie algebras into
an example of hom-Lie algebras. In Section \ref{sec:extofhlas}, for the class of hom-Lie algebras, we develop
the central extension theory, providing homological type conditions useful for showing existence of central
extensions and for their construction. Here, we required that the central extension of a hom-Lie algebra is also
a hom-Lie algebra. In particular, the standard theory of central extensions of Lie algebras becomes a natural
special case of the theory for hom-Lie algebras when no non-identity twisting is present. In particular, this
implies that in the specific examples of deformation families of Witt and Virasoro type algebras constructed
within this framework, the corresponding non-deformed Witt and Virasoro type Lie algebras are included as the
algebras corresponding to those specific values of deformation parameters which remove the non-trivial twisting.
In Section \ref{sec:cextqVirqWittlinear}, we demonstrate the use of the central extension theory for hom-Lie
algebras by applying it to the construction of a central hom-Lie algebra extension of the $q$-deformed Witt
algebra from Theorem \ref{th:qWittlinear}, which is a $q$-deformation of Virasoro Lie algebra. For $q=1$ one
indeed recovers the usual Virasoro Lie algebra as is expected from our general approach.

It should be mentioned that the use of $q$-deformed Jacobi identities for constructing $q$-deformations of the
Witt and Virasoro algebras has been considered in physical and mathematical literature before
\cite{AizSato1,ChaiIsLukPopPres,ChaiKuLuk,ChakJag,Polychronakos,SatoHT1,SatoHT2,SatoHT3,SatoHT4}. Of these
works, the closest to our results on hom-Lie algebras comes \cite{AizSato1} where the two identities,
skew-symmetry and a twisted from inside three-term Jacobi identity, almost as the one for hom-Lie algebras, have
been clearly stated as a definition of a class of non-associative algebras, and then used as the conditions
required to be satisfied by the central extension of a \mbox{$q$-deformation} of the Witt algebra from
\cite{CurtrZachos1}. This results in a $q$-deformation of the Virasoro Lie algebra somehow related to that in
the example we described in section \ref{sec:cextqVirqWittlinear}. Whether a particular deformation of the Witt
or Virasoro algebra obtained by various constructions satisfy some kinds of Jacobi type identities is considered
to be an important problem. The generalized twisted $6$-term Jacobi identity obtained in our construction, gives
automatically by specialization the deformed Jacobi identities satisfied by the corresponding particular
deformations of the Witt and Virasoro algebras. There are also works employing usual and super Jacobi identities
as conditions on central extensions and their deformations (for example
\cite{AvanFrappatRossiSorba,FairNuyZach,KemmokuSatoHT,ZhaCZZhaoWZ}). Putting these works within context of our
approach would be of interest.

We would also like to note that in the works \cite{BlochS1,Li,Su1,Su2}, in the case of usual derivations on
Laurent polynomials, it has been specifically noted that a Lie bracket can be defined by expressions somewhat
resembling a special case of \eqref{eq:GenWittProdFormula}. We also would like to mention that $q$-deformations
of the Witt and Virasoro algebras were considered indirectly as an algebra of pseudo \mbox{$q$-difference}
operators based on the $q$-derivative on Laurent polynomials in \cite{KacRadul,Kassel1,KhesinLyubRoger}. We
believe that it should be possible, and would be of direct interest, to extend the results of these works to our
general context of \mbox{$\sigma$-derivations,} and we hope to contribute to this cause in future work. For the
reader's convenience, we  have also included in the bibliography, without further reference in the text, some
works we know of, concerned with other specific examples of deformations of Witt algebras that we believe could
be considered in our framework, leaving the possibility of this as an open question for the moment.

We also feel that the further development should include using our construction for building more examples of
deformed or twisted Witt and Virasoro type algebras based on differential and difference type operators on
function spaces studied extensively in analysis and in numerical mathematics, and on functions on algebraic
varieties important in algebraic geometry and its applications. It could be of interest to extend our
constructions and examples over fields of finite characteristic, or various number fields. Development of the
representation theory for the parametric families of Witt and Virasoro type algebras arising within our method,
and understanding to which extent the representations of non-deformed Witt and Virasoro algebras appear as limit
points will be important for applications in physics.

\section{Some general considerations}

\subsection{Generalized derivations on commutative algebras and on UFD's}

We begin with some definitions. Throughout this section, $\A$ is
an associative $\C$-algebra, and $\s$ and $\t$ denote two
different algebra endomorphisms on $\A$.
\begin{dfn}
 A \textit{$(\s,\t)$-derivation} $D$ on $\A$ is a
    $\C$-linear map satisfying
    $$D(ab)=D(a)\t(b)+\s(a)D(b),$$ where $a,b\in\A$. The set of all
$(\s,\t)$-derivations on $\A$ is denoted by
    $\mfD_{(\s,\t)}(\A)$.
\end{dfn}
\begin{dfn}
    With notation as above, a \textit{$\s$-derivation} on $\A$ is a
$(\s,\id)$-derivation,
    i.e. a $\C$-linear map $D$ satisfying
    $$D(ab)=D(a)b+\s(a)D(b),$$ for $a,b\in\A$. We denote the set of all
$\s$-derivations by $\mfD_\s(\A)$.
\end{dfn}
 From now on, when speaking of unique factorization domains (UFD), we shall always
mean a commutative associative algebra over $\C$ with unity $1$
and with no zero-divisors, such that any element can be written in
a unique way (up to a multiple of an invertible element) as a
product of irreducible elements, i.e. elements which cannot be
written as a product of two non-invertible elements. Examples of
unique factorization domains include $\C[x_1,\ldots,x_n]$, and the
algebra $\C[t,t^{-1}]$ of Laurent polynomials.

When $\s(x)a=a\s(x)$ (or $\t(x)a=a\t(x)$) for all $x,a\in \A$ and
in particular when $\A$ is commutative,
$\mathfrak{D}_{(\s,\t)}(\A)$ carries a natural left (or right)
$\A$-module structure by $(a,D)\mapsto a\cdot D: x \mapsto aD(x)$.
If $a,b\in\A$ we shall write $a\big|b$ if there is an element
$c\in\A$ such that $ac=b$. If $S\subseteq\A$ is a subset of $\A$,
a greatest common divisor, $\GCD(S)$, of $S$ is defined as an
element of $\A$ satisfying
\begin{align}\label{eqn:polufd3}
    \GCD(S) \big| a \quad \text{for } a\in S,
\end{align} and
\begin{align} \label{eqn:polufd4}
    b\big| a \quad \text{for } a\in S \Longrightarrow b\big| \GCD(S).
\end{align}
It follows directly from the definition that
\begin{align}\label{eq:GCDsubset}
    S\subseteq T\subseteq\A \quad\Longrightarrow\quad \GCD(T)\big|\GCD(S)
\end{align}whenever $\gcd(S)$ and $\gcd(T)$ exist.
If $\A$ is a unique factorization domain one can show that a
$\GCD(S)$ exists for any nonempty subset $S$ of $\A$ and that this
element is unique up to a multiple of an invertible element in
$\A$. Thus we are allowed to speak of \emph{the} $\gcd$.

\begin{lem} \label{lem:pollem3}
Let $\A$ be a commutative algebra. Let $\s$ and $\t$ be two
algebra endomorphisms on $\A$, and let $D$ be a
$(\s,\t)$-derivation on $\A$. Then
$$D(x)(\t(y)-\s(y))=0$$ for all $x\in\ker(\t-\s)$ and $y\in\A$. Moreover, if $\A$
has no zero-divisors and $\s\neq\t$, then
\begin{align} \label{eqn:pollem3eq}
    \ker (\t-\s)\subseteq \ker D.
\end{align}
\end{lem}
\begin{proof}
Let $y\in \A$ and let $x\in\ker(\t-\s)$. Then
\begin{align*}
    0&=D(xy-yx)=D(x)\t(y)+\s(x)D(y)-D(y)\t(x)-\s(y)D(x)=\\
    &=D(x)(\t(y)-\s(y))-D(y)(\t(x)-\s(x))=D(x)(\t(y)-\s(y)).
\end{align*}Furthermore, if $\A$ has no zero-divisors and if there is a $y\in\A$
such that $\t(y)\neq\s(y)$ then $D(x)=0$.
\end{proof}

\begin{thm} \label{thm:polndifferent}
Let $\s$ and $\t$ be different algebra endomorphisms on a unique
factorization domain $\A$. Then $\mathfrak{D}_{(\s,\t)}(\A)$ is
free of rank one as an $\A$-module with generator
\begin{align} \label{eqn:polthmn}
    \De:=\frac{(\t-\s)}{g}\quad:\quad x\longmapsto \frac{(\t-\s)(x)}{g}
\end{align} where $g=\GCD\big((\t-\s)(\A)\big)$.
\end{thm}
\begin{proof}
We note first that $(\t-\s)/g$ is a $(\s,\t)$-derivation on $\A$:
\begin{align*}
    \frac{(\t-\s)(xy)}{g}&=\frac{\t(x)\t(y)-\s(x)\s(y)}{g}=\\
    &=\frac{\big(\t(x)-\s(x)\big)\t(y)+\s(x)\big(\t(y)-\s(y)\big)}{g}=\\
    &=\frac{(\t-\s)(x)}{g}\cdot\t(y)+\s(x)\cdot\frac{(\t-\s)(y)}{g},
\end{align*}
for $x,y\in\A$. Next we show that $(\t-\s)/g$ generates a free
$\A$-module of rank one. So suppose that
\begin{align}\label{eq:polthmFree1}
    x\cdot\frac{\t-\s}{g}=0,
\end{align} for some $x\in\A$. Since $\t\neq\s$, there is an $y\in\A$ such that
$(\t-\s)(y)\neq 0$. Application of both sides in
(\ref{eq:polthmFree1}) to this $y$ yields
$$x\cdot\frac{(\t-\s)(y)}{g}=0.$$
Since $\A$ has no zero-divisors, it then follows that $x=0$. Thus
$$\A\cdot\frac{\t-\s}{g}$$ is a free $\A$-module of rank one.

It remains to show that
$\mathfrak{D}_{(\s,\t)}(\A)\subseteq\A\cdot\frac{\t-\s}{g}$. Let
$D$ be a $(\s,\t)$-derivation on $\A$. We want to find $a_D\in\A$
such that
\begin{align}\label{eq:polthmGREJ}
    D(x)=a_D\cdot\frac{(\t-\s)(x)}{g}
\end{align}for $x\in\A$. We will define
\begin{align} \label{eqn:polthmwant}
    a_D=\frac{D(x)\cdot g}{(\t-\s)(x)}
\end{align} for some $x$ such that $(\t-\s)(x)\neq 0$. For this to be
possible, we must show two things. First of all, that
\begin{align}\label{eqn:polthmmust1}
    (\t-\s)(x) \;\big|\;D(x)\cdot g \quad\text{for any } x \text{ with }
(\t-\s)(x)\neq 0
\end{align} and secondly, that
\begin{align}\label{eqn:polthmmust2}
    \frac{D(x)\cdot g}{(\t-\s)(x)}=\frac{D(y)\cdot g}{(\t-\s)(y)} \quad\text{for any }
    x,y \text{ with } (\t-\s)(x), (\t-\s)(y)\neq 0.
\end{align} Suppose for a moment that (\ref{eqn:polthmmust1}) and
(\ref{eqn:polthmmust2}) were true. Then it is clear that if we
define $a_D$ by (\ref{eqn:polthmwant}), the formula
(\ref{eq:polthmGREJ}) holds for any $x\in\A$ satisfying
$(\t-\s)(x)\neq 0$. But (\ref{eq:polthmGREJ}) also holds when
$x\in\A$ is such that $(\t-\s)(x)=0$, because then $D(x)=0$ also,
by Lemma \ref{lem:pollem3}.

We first prove (\ref{eqn:polthmmust1}). Let $x,y\in\A$ be such
that $(\t-\s)(x), (\t-\s)(y)\neq 0$. Then we have
\begin{align*}
    0&=D(xy-yx)=D(x)\t(y)+\s(x)D(y)-D(y)\t(x)-\s(y)D(x)=\\
    &=D(x)(\t(y)-\s(y))-D(y)(\t(x)-\s(x)),
\end{align*}
so that
\begin{align} \label{eqn:polthmdiv}
    D(x)(\t(y)-\s(y))=D(y)(\t(x)-\s(x)).
\end{align} Now define a function
$h:\A\times\A\to\A$ by setting
$$h(z,w)=\GCD(\t(z)-\s(z),\t(w)-\s(w))\qquad\text{for } z,w\in\A.$$
By the choice of $x$ and $y$, we have $h(x,y)\neq 0$. Divide both
sides of (\ref{eqn:polthmdiv}) by $h(x,y)$:
\begin{align} \label{eqn:polthmdiv2}
    D(x)\frac{\t(y)-\s(y)}{h(x,y)}= D(y)\frac{\t(x)-\s(x)}{h(x,y)}.
\end{align}
By construction 
$$\GCD\Big (\frac{\t(y)-\s(y)}{h(x,y)},\frac{\t(x)-\s(x)}{h(x,y)}\Big )=1.$$
Therefore, using that $\A$ is a UFD, we deduce from
(\ref{eqn:polthmdiv2}) that
$$\frac{\t(x)-\s(x)}{h(x,y)} \;\big|\; D(x),$$
implying that
\begin{align} \label{eqn:polthmdiv3}
    (\t-\s)(x) \;\big|\; D(x)\cdot h(x,y)
\end{align} for any $x,y\in\A$ with $(\t-\s)(x), (\t-\s)(y)\neq 0$. Let
$S=\A\backslash\ker(\t-\s)$. Then from (\ref{eqn:polthmdiv3}) and
property (\ref{eqn:polufd4}) of the $\GCD$ we get
\begin{align} \label{eqn:polthmdiv5}
    (\t-\s)(x) \;\big|\;   D(x)\cdot \GCD(h(x,S))
\end{align} for all $x\in\A$ with $(\t-\s)(x)\neq 0$. But
\begin{align*}
    \GCD(h(x,S))&=\GCD\Big(\big\{\GCD\big((\t-\s)(x),(\t-\s)(s)\big) \; |\; s\in
S\big\}\Big)=\\
    &=\GCD\big((\t-\s)(S)\cup\{(\t-\s)(x)\}\big)=\\
    &=\GCD\big((\t-\s)(\A)\cup\{(\t-\s)(x)\}\big)=\\
    &=g.
\end{align*}
Thus (\ref{eqn:polthmdiv5}) is equivalent to
(\ref{eqn:polthmmust1}) which was to be proved.

Finally, we prove (\ref{eqn:polthmmust2}). Let $x,y\in\A$ be such
that $(\t-\s)(x), (\t-\s)(y)\neq 0$. Then
\begin{align*}
    0&=D(xy-yx)=D(x)\t(y)+\s(x)D(y)-D(y)\t(x)-\s(y)D(x)=\\
    &=D(x)(\t(y)-\s(y))-D(y)(\t(x)-\s(x)),
\end{align*}
which, after multiplication by $g$ and division by
$(\t-\s)(x)\cdot(\t-\s)(y)$ proves (\ref{eqn:polthmmust2}). This
completes the proof of the existence of $a_D$, and hence the proof
of the theorem.
\end{proof}

\subsection{A bracket on $\s$-derivations}
\label{sec:generalizationOfWitt} The Witt algebra is isomorphic to
the Lie algebra $\mathfrak{D}(\C[t,t^{-1}])$ of all derivations of
the commutative unital algebra of all complex Laurent polynomials:
$$\C[t,t^{-1}]=\{\sum_{k\in\Z}a_k t^k\;|\; a_k\in\C, \text{only finitely many
non-zero}\}.$$ In this section we will use this fact as a starting
point for a generalization of the Witt algebra to an algebra
consisting of $\s$-derivations.

We let $\A$ be a commutative associative algebra over $\C$ with
unity $1$, as in the example $\A=\C[t,t^{-1}]$ from the previous
paragraph. When we speak of homomorphisms (endomorphisms) in the
sequel we will always mean \emph{algebra} homomorphisms
(endomorphisms), except where otherwise indicated.
 If \mbox{$\s:\A\to\A$} is a homomorphism of algebras, we denote, as before, the
\mbox{$\A$-module} of all $\s$-derivations on $\A$
 by $\mfD_\s(\A)$. For clarity we will denote the module multiplication by $\cdot$
and the algebra multiplication in $\A$ by juxtaposition. The
\emph{annihilator} $\Ann(D)$ of an element $D\in\mfD_\s(\A)$ is
the set of all $a\in\A$ such that $a\cdot D=0$. It is easy to see
that $\Ann(D)$ is an ideal in $\A$ for any $D\in\mfD_\s(\A)$.

We now fix a homomorphism $\s:\A\to\A$, an element
$\De\in\mfD_\s(\A)$, and an element $\d\in\A$, and we assume that
these objects satisfy the following two conditions:
\begin{align}\label{eq:GenWittCond1}
    \s(\Ann(\De))\subseteq \Ann(\De),
\end{align}
\begin{align}\label{eq:GenWittCond2}
    \De(\s(a)) = \d\s(\De(a)),\quad\text{for }a\in\A.
\end{align}
Let
$$\A\cdot\De=\{a\cdot\De\;|\;a\in\A\}$$
denote the cyclic $\A$-submodule of $\mfD_\s(\A)$ generated by
$\De$. We have the following theorem, which introduces a
$\C$-algebra structure on $\A\cdot\De$.
\begin{thm} \label{thm:GenWitt}
If (\ref{eq:GenWittCond1}) holds then the map
\begin{align*}
    [\cdot,\cdot]_\s:\A\cdot\De\times\A\cdot\De \to
    \A\cdot\De
\end{align*} defined by setting
\begin{align} \label{eq:GenWittProdDef}
    [a\cdot\De,b\cdot\De]_\s=(\s(a)\cdot\De)\circ(b\cdot\De)-(\s(b)\cdot\De)
    \circ(a\cdot\De), \quad\text{for }a,b\in\A,
\end{align}
where $\circ$ denotes composition of functions, is a well-defined
$\C$-algebra product on the $\C$-linear space $\A\cdot\De$, and it
satisfies the following identities for $a,b,c\in\A$:
\begin{align}\label{eq:GenWittProdFormula}
    [a\cdot\De, b\cdot\De]_\s=\big(\s(a)\De(b)-\s(b)\De(a)\big)\cdot\De,
\end{align}
\begin{align}\label{eq:GenWittSkew}
    [a\cdot\De, b\cdot\De]_\s=-[b\cdot\De, a\cdot\De]_\s.
\end{align} In addition, if (\ref{eq:GenWittCond2}) holds, then
\begin{equation} \label{eq:GenWittJacobi}
\begin{split}
    [\s&(a)\cdot\De,[b\cdot\De,c\cdot\De]_\s]_\s+\d\cdot[a\cdot\De,[b\cdot\De,c\cdot\De]_\s]_\s+\\
    +&[\s(b)\cdot\De,[c\cdot\De,a\cdot\De]_\s]_\s+\d\cdot[b\cdot\De,[c\cdot\De,a\cdot\De]_\s]_\s+\\
    +&[\s(c)\cdot\De,[a\cdot\De,b\cdot\De]_\s]_\s+\d\cdot[c\cdot\De,[a\cdot\De,b\cdot\De]_\s]_\s=0.
\end{split}
\end{equation}
\end{thm}
\begin{remark}\label{rem1}
    An important thing to notice is that the bracket $[\cdot,\cdot]_\s$ defined in
the theorem
    depends on the generator $\De$ of
    the cyclic submodule $\A\cdot\De$ of $\mfD_\s(\A)$ in an essential way. This
reveals that one should in fact write
    $[\cdot,\cdot]_{\s,\De}$ to explicitly indicate which $\De$ is chosen.
    Suppose, however, we choose another generator $\De'$ of $\A\cdot\De$. Then
$\De'=u\De$ for an element $u\in\A$
    (not necessarily a unit). Take
    elements $a\cdot\De'\,,b\cdot\De'\in\A\cdot\De$. Then the following calculation
shows how two different
    brackets relate when changing the generator (we use the commutativity of $\A$
freely):
    \begin{multline*}
        \s(u)[a\cdot\De', b\cdot\De']_{\s,\De'}=\text{[The definition of the
bracket]}=\\
        =\big (\s(a) u\s(u)\cdot\De\big )\circ (bu\De)-\big
(\s(b)u\s(u)\cdot\De)\circ (au\cdot\De)=\\
        =u\cdot\Big (\big (\s(au)\cdot\De\big )\circ\big(bu\cdot\De)-\big
(\s(bu)\cdot\De\big )\circ\big(au\cdot\De)
        \Big)=\\
        =u\cdot[au\cdot\De,bu\cdot\De]_{\s,\De}=u\cdot[a\cdot\De',b\cdot\De']_{\s,\De}
    \end{multline*}so the ''base change''-relation is
    $$\s(u)\cdot[a\cdot\De',
b\cdot\De']_{\s,\De'}=u\cdot[a\cdot\De',b\cdot\De']_{\s,\De}.$$
For the most part
    of this paper, we have a fixed
    generator and so we suppress the dependence on the generator from the bracket
notation and simply write
    $[\cdot,\cdot]_{\s}.$
On the other hand, if $\A$ has no zero-divisors, we shall see
later in Proposition \ref{prop72} that the dependence of the
generator $\De$ is not essential.
\end{remark}
\begin{remark}
    The identity (\ref{eq:GenWittProdFormula}) is just a formula expressing the
product defined in
    (\ref{eq:GenWittProdDef}) as an element of $\A\cdot\De$. Identities
(\ref{eq:GenWittSkew}) and
    (\ref{eq:GenWittJacobi}) are more essential, expressing, respectively,
skew-symmetry and a generalized ($(\s,\d)$-twisted) Jacobi identity for the product defined by
(\ref{eq:GenWittProdDef}).
\end{remark}
Before coming to the proof of the theorem we introduce a
convenient notation. If $f:\A\times\A\times\A\to \A\cdot\De$ is a
function, we will write
$$\cs_{a,b,c}f(a,b,c)$$
for the cyclic sum
$$f(a,b,c)+f(b,c,a)+f(c,a,b).$$
We note the following properties of the cyclic sum:
$$\cs_{a,b,c}\big(x\cdot f(a,b,c)+y\cdot g(a,b,c)\big)=x\cdot\cs_{a,b,c}f(a,b,c)+
y\cdot\cs_{a,b,c}g(a,b,c),$$
$$\cs_{a,b,c}f(a,b,c)=\cs_{a,b,c}f(b,c,a)=
\cs_{a,b,c}f(c,a,b), $$ where
$f,g:\A\times\A\times\A\to\A\cdot\De$ are two functions, and
$x,y\in\A$. Combining these two identities we obtain
\begin{align}\label{shift_prop}
    \cs_{a,b,c}\big(f(a,b,c)+g(a,b,c)\big)&=
    \cs_{a,b,c}\big(f(a,b,c)+g(b,c,a)\big)=\notag\\
    &= \cs_{a,b,c}\big(f(a,b,c)+g(c,a,b)\big).
\end{align} With this notation, (\ref{eq:GenWittJacobi}) can be written
\begin{align}\label{eq:GenWittJacobiCycNot}
    \cs_{a,b,c}\Big\{[\s(a)\cdot\De,[b\cdot\De,c\cdot\De]_\s]_\s
  +\d\cdot [a\cdot\De,[b\cdot\De,c\cdot\De]_\s]_\s\Big\}=0.
\end{align}
We now turn to the proof of Theorem \ref{thm:GenWitt}.
\begin{proof}
We must first show that $[\cdot,\cdot]_\s$ is a well defined
function. That is, if $a_1\cdot\De=a_2\cdot\De$, then
\begin{align} \label{eq:GenWittPropPrf1}
    [a_1\cdot\De, b\cdot\De]_\s=[a_2\cdot\De, b\cdot\De]_\s,
\end{align} and
\begin{align} \label{eq:GenWittPropPrf2}
    [b\cdot\De, a_1\cdot\De]_\s=[b\cdot\De, a_2\cdot\De]_\s,
\end{align} for $b\in\A$. Now
$a_1\cdot\De=a_2\cdot\De$ is equivalent to $a_1-a_2\in\Ann(\De)$.
Therefore, using the assumption (\ref{eq:GenWittCond1}), we also
have $\s(a_1-a_2)\in\Ann(\De)$. Hence
\begin{multline*}
    [a_1\cdot\De, b\cdot\De]_\s-[a_2\cdot\De, b\cdot\De]_\s =
    (\s(a_1)\cdot\De)\circ(b\cdot\De)-(\s(b)\cdot\De)\circ(a_1\cdot\De)-\\
    \quad -(\s(a_2)\cdot\De)\circ(b\cdot\De)+(\s(b)\cdot\De)\circ(a_2\cdot\De)=\\
    =(\s(a_1-a_2)\cdot\De)\circ(b\cdot\De)-(\s(b)\cdot\De)\circ((a_1-a_2)\cdot\De)=0,
\end{multline*}
which shows (\ref{eq:GenWittPropPrf1}). The proof of
(\ref{eq:GenWittPropPrf2}) is analogous.

Next we prove (\ref{eq:GenWittProdFormula}), which also shows that
$\A\cdot\De$ is closed under $[\cdot,\cdot]_\s$. Let $a,b,c\in\A$
be arbitrary. Then, since $\De$ is a $\s$-derivation on $\A$ we
have
\begin{multline*}
    [a\cdot\De, b\cdot\De]_\s(c)=
    (\s(a)\cdot\De)\Big((b\cdot\De)(c)\Big)-(\s(b)\cdot\De)\Big((a\cdot\De)(c)\Big)=\\
    =\s(a)\De\big(b\De(c)\big)-\s(b)\De\big(a\De(c)\big)=\\
    =\s(a)\big(\De(b)\De(c)+\s(b)\De(\De(c))\big)-\s(b)\big(\De(a)\De(c)+\s(a)\De(\De(c))\big)=\\
    =\big(\s(a)\De(b)-\s(b)\De(a)\big)\De(c)+(\s(a)\s(b)-\s(b)\s(a))\De(\De(c)).
\end{multline*}
Since $\A$ is commutative, the last term is zero. Thus
(\ref{eq:GenWittProdFormula}) is true. The skew-symmetry identity
(\ref{eq:GenWittSkew}) is clear from the definition
(\ref{eq:GenWittProdDef}). Using the linearity of $\s$ and $\De$,
and the definition of $[\cdot,\cdot]_\s$, or the formula
(\ref{eq:GenWittProdFormula}), it is also easy to see that
$[\cdot,\cdot]_\s$ is bilinear.

It remains to prove (\ref{eq:GenWittJacobi}). Using
(\ref{eq:GenWittProdFormula}) and that $\De$ is a $\s$-derivation
on $\A$ we get
\begin{align}
    [&\s(a)\cdot\De,[b\cdot\De,c\cdot\De]_\s]_\s=[\s(a)\cdot\De,\big(\s(b)\De(c)-\s(c)\De(b)\big)\cdot\De]_\s=\nonumber\\
    &=\Big\{\s^2(a)\De\big(\s(b)\De(c)-\s(c)\De(b)\big)-\nonumber\\
    &\quad-\s\big(\s(b)\De(c)-\s(c)\De(b)\big)\De\big(\s(a)\big)\Big\}\cdot\De=\nonumber\\
    &=\Big\{\s^2(a)\Big(\De\big(\s(b)\big)\De(c)+\s^2(b)\De^2(c)-
    \De\big(\s(c)\big)\De(b)-\s^2(c)\De^2(b)\Big)-\nonumber\\
    &\quad-\Big(\s^2(b)\s\big(\De(c)\big)-\s^2(c)\s\big(\De(b)\big)\Big)
    \De\big(\s(a)\big)\Big\}\cdot\De=\nonumber\\
    &=\s^2(a)\De\big(\s(b)\big)\De(c)\cdot\De+\s^2(a)\s^2(b)\De^2(c)\cdot\De-\nonumber\\
    &\quad-\s^2(a)\De\big(\s(c)\big)\De(b)\cdot\De-\s^2(a)\s^2(c)\De^2(b)\cdot\De-\nonumber\\
    &\quad-\s^2(b)\s\big(\De(c)\big)\De\big(\s(a)\big)\cdot\De
    +\s^2(c)\s\big(\De(b)\big)\De\big(\s(a)\big)\cdot\De,\label{eq:GenWittCalc1}
\end{align}
where $\s^2=\s\circ\s$ and $\De^2=\De\circ\De$. Applying cyclic
summation to the second and fourth term in (\ref{eq:GenWittCalc1})
we get
\begin{multline*}
    \cs_{a,b,c}\Big\{\s^2(a)\s^2(b)\De^2(c)\cdot\De-\s^2(a)\s^2(c)\De^2(b)\cdot\De\Big\}=\\
    =\cs_{a,b,c}\Big\{\s^2(a)\s^2(b)\De^2(c)\cdot\De-\s^2(b)\s^2(a)\De^2(c)\cdot\De\Big\}=0,
\end{multline*}
using (\ref{shift_prop}) and that $\A$ is commutative. Similarly,
if we apply cyclic summation to the fifth and sixth term in
(\ref{eq:GenWittCalc1}) and use the relation
(\ref{eq:GenWittCond2}) we obtain
\begin{multline*}
    \cs_{a,b,c}\Big\{-\s^2(b)\s\big(\De(c)\big)\De\big(\s(a)\big)\cdot\De
     +\s^2(c)\s\big(\De(b)\big)\De\big(\s(a)\big)\cdot\De\Big\}=\\
    =\cs_{a,b,c}\Big\{-\s^2(b)\s\big(\De(c)\big)\d\s\big(\De(a)\big)\cdot\De
     +\s^2(c)\s\big(\De(b)\big)\d\s\big(\De(a)\big)\cdot\De\Big\}=\\
    =\d\cdot\cs_{a,b,c}\Big\{-\s^2(b)\s\big(\De(c)\big)\s\big(\De(a)\big)\cdot\De
     +\s^2(b)\s\big(\De(a)\big)\s\big(\De(c)\big)\cdot\De\Big\}=0,
\end{multline*}
where we again used (\ref{shift_prop}) and the commutativity of
$\A$. Consequently, the only terms in the right hand side of
(\ref{eq:GenWittCalc1}) which do not vanish when we take cyclic
summation are the first and the third. In other words,
\begin{multline}\label{eq:GenWittSteg}
    \cs_{a,b,c}[\s(a)\cdot\De,[b\cdot\De,c\cdot\De]_\s]_\s=\\
    =\cs_{a,b,c}\Big\{\s^2(a)\De\big(\s(b)\big)\De(c)\cdot\De-
     \s^2(a)\De\big(\s(c))\De(b)\cdot\De\Big\}.\quad
\end{multline}

We now consider the other term in (\ref{eq:GenWittJacobiCycNot}).
First note that from (\ref{eq:GenWittProdFormula}) we have
$$[b\cdot\De,c\cdot\De]_\s=\big(\De(c)\s(b)-\De(b)\s(c)\big)\cdot\De$$
since $\A$ is commutative. Using first this and then
(\ref{eq:GenWittProdFormula}) we get
\begin{align*}
    \d&\cdot[a\cdot\De,[b\cdot\De,c\cdot\De]_\s]_\s=\\
    &=\d\cdot[a\cdot\De,\big(\De(c)\s(b)-\De(b)\s(c)\big)\cdot\De]_\s=\\
    &=\d\Big(\s(a)\De\big(\De(c)\s(b)-\De(b)\s(c)\big)
    -\s\big(\De(c)\s(b)-\De(b)\s(c)\big)\De(a)\Big)\cdot\De=\\
    &=\d\Big\{\s(a)\Big(\De^2(c)\s(b)+\s\big(\De(c)\big)\De\big(\s(b)\big)
    -\De^2(b)\s(c)-\s\big(\De(b)\big)\De\big(\s(c)\big)\Big)\\
    &\quad-\Big(\s\big(\De(c)\big)\s^2(b)-\s\big(\De(b)\big)\s^2(c)\Big)\De(a)\Big\}\cdot\De=\\
    &=\d\s(a)\De^2(c)\s(b)\cdot\De+\d\s(a)\s\big(\De(c)\big)\De\big(\s(b)\big)\cdot\De-\\
    &\quad-\d\s(a)\De^2(b)\s(c)\cdot\De-\d\s(a)\s\big(\De(b)\big)\De\big(\s(c)\big)\cdot\De-\\
    &\quad-\d\s\big(\De(c)\big)\s^2(b)\De(a)\cdot\De+\d\s\big(\De(b)\big)\s^2(c)\De(a)\cdot\De.
\end{align*}
Using (\ref{eq:GenWittCond2}), this is equal to
\begin{multline*}
    \quad\d\s(a)\De^2(c)\s(b)\cdot\De+\s(a)\De\big(\s(c)\big)\De\big(\s(b)\big)\cdot\De-\\
    \quad-\d\s(a)\De^2(b)\s(c)\cdot\De-\s(a)\De\big(\s(b)\big)\De\big(\s(c)\big)\cdot\De-\\
    \quad-\De\big(\s(c)\big)\s^2(b)\De(a)\cdot\De+\De\big(\s(b)\big)\s^2(c)\De(a)\cdot\De=\\
    =\d\s(a)\De^2(c)\s(b)\cdot\De-\d\s(a)\De^2(b)\s(c)\cdot\De-\\
    \quad-\De\big(\s(c)\big)\s^2(b)\De(a)\cdot\De+\De\big(\s(b)\big)\s^2(c)\De(a)\cdot\De.
\end{multline*}
The first two terms of this last expression vanish after a cyclic
summation and using (\ref{shift_prop}), so we get
\begin{multline}\label{eq:GenWittSteg2}
    \cs_{a,b,c}\d\cdot[a\cdot\De,[b\cdot\De,c\cdot\De]_\s]_\s=\\
    =\cs_{a,b,c}\Big\{-
    \De\big(\s(c)\big)\s^2(b)\De(a)\cdot\De+\De\big(\s(b)\big)\s^2(c)\De(a)\cdot\De\Big\}.\quad
\end{multline}
Finally, combining this with (\ref{eq:GenWittSteg}) we deduce
\begin{align*}
    \cs_{a,b,c}&\Big\{[\s(a)\cdot\De,[b\cdot\De,c\cdot\De]_\s]_\s
    +\d[a\cdot\De,[b\cdot\De,c\cdot\De]_\s]_\s\Big\}=\\
    &=\cs_{a,b,c}[\s(a)\cdot\De,[b\cdot\De,c\cdot\De]_\s]_\s
    +\cs_{a,b,c}\d[a\cdot\De,[b\cdot\De,c\cdot\De]_\s]_\s=\\
    &=\cs_{a,b,c}\Big\{\s^2(a)\De\big(\s(b)\big)\De(c)\cdot\De-
    \s^2(a)\De\big(\s(c))\De(b)\cdot\De\Big\}+\\
    &\quad+\cs_{a,b,c}\Big\{-\De\big(\s(c)\big)\s^2(b)\De(a)\cdot\De
    +\De\big(\s(b)\big)\s^2(c)\De(a)\cdot\De\Big\}=\\
    &=\cs_{a,b,c}\Big\{\s^2(a)\De\big(\s(b)\big)\De(c)\cdot\De-
    \s^2(a)\De\big(\s(c))\De(b)\cdot\De\Big\}+\\
    &\quad+\cs_{a,b,c}\Big\{-\De\big(\s(b)\big)\s^2(a)\De(c)\cdot\De
    +\De\big(\s(c)\big)\s^2(a)\De(b)\cdot\De\Big\}=\\
    &=0,
\end{align*}
as was to be shown. The proof is complete.
\end{proof}
\begin{remark}
If $\A$ is not assumed to be commutative, the construction still
works if one impose on $\De$ the additional condition that
$$[a,b]\De(c)=0\quad\text{for all } a,b,c\in\A.$$
Then the mapping $x\cdot\De:b\mapsto x\De(b)$ is again a
$\s$-derivation for all $x\in\A$. As before $\A\cdot\De$ is a left
$\A$-module. Then Theorem \ref{thm:GenWitt} remain valid with the
same proof. We only need to note that, although $\A$ is not
commutative we have $[a,b]\cdot\De=0$ which is to say that
$$ab\cdot\De=ba\cdot\De.$$
\end{remark}
\begin{prop}\label{prop72}
    If $\A$ is a commutative $\C$-algebra without zero-divisors, and if $0\neq
\De\in\mfD_\s(\A)$ and $0\neq \De'\in\mfD_\s(\A)$
     generates the same cyclic $\A$-submodule $\M$ of $\mfD_\s(\A)$, where
$\s:\A\longrightarrow \A$ is an algebra
     endomorphism, then there is a unit $u\in\A$ such that
     \begin{align}\label{eq:uDD'}
        [x,y]_{\s,\De}=u[x,y]_{\s,\De'}.
     \end{align}Furthermore, if $u\in\C$ then
     \begin{align*}(\M,[\,\cdot,\,\cdot]_{\s,\De})\cong
(\M,[\,\cdot,\,\cdot]_{\s,\De'})
     \end{align*}
\end{prop}
\begin{proof}That $\De$ and $\De'$ generates the same cyclic submodule implies that
there are $u_1, u_2$ such that $\De=u_1\De'$ and $\De'=u_2\De$.
This means that $u_1u_2\De=u_1\De'=\De$ or equivalently
$(u_1u_2-1)\De=0$. Choose $a\in\A$ such that $\De(a)\neq 0$. Then
$(u_1u_2-1)\De(a)=0$ implies that $u_1u_2-1=0$ and so $u_1$ and
$u_2$ are both units. We now use Remark \ref{rem1} to get
$$\s(u_2)\cdot [x,y]_{\s,\De'}=u_2\cdot [x,y]_{\s,\De}.$$ Then $u=\s(u_2)/u_2$
satisfies (\ref{eq:uDD'}). Now, if $u\in\C$ define $\varphi :
(\M,[\,\cdot,\,\cdot]_{\s,\De})\longrightarrow
(\M,[\,\cdot,\,\cdot]_{\s,\De'})$ by $\varphi(x)=ux$. Then
$$\varphi([x,y]_{\s,\De})=u[x,y]_{\s,\De}=u^2[x,y]_{\s,\De'}=[ux,uy]_{\s,\De'}=[\varphi(x),\varphi(y)]_{\s,\De'}.$$
\end{proof}
\begin{dfn} \label{def:genWittalg}
    Let $\A$ be commutative and associative algebra, $\s:\A\to\A$ an algebra
endomorphism and $\De$ a
    $\s$-derivation on $\A$. Then, a \textit{$(\A,\s,\De)$-Witt algebra} (or a
\textit{generalized Witt
    algebra}) is the non-associative algebra $(\A \cdot\De,[\cdot,\cdot]_{\s,\De})$
with the product defined by
    $$[a\cdot\De,b\cdot\De]_{\s,\De}=\big (\s(a)\De(b)-\s(b)\De(a)\big )\cdot\De.$$
\end{dfn}
\begin{example}
Take $\A=\C[t,t^{-1}]$, $\s=\id_\A$, the identity operator on
$\A$, $\De=\frac{d}{dt}$, and $\d=1$. In this case one can show
that $\A\cdot\De$ is equal to the whole $\mathfrak{D}_\s(\A)$. The
conditions (\ref{eq:GenWittCond1}) and (\ref{eq:GenWittCond2}) are
trivial to check. The definition (\ref{eq:GenWittProdDef})
coincides with the usual Lie bracket of derivations, and equation
(\ref{eq:GenWittJacobi}) reduces to twice the usual Jacobi
identity. Hence we recover the ordinary Witt algebra.
\end{example}
\begin{example}\label{ex:delta} Let $\A$ be a unique factorization domain, and let
$\s:\A\to\A$ be a homomorphism, different from the identity. Then
by Theorem \ref{thm:polndifferent},
$$\mathfrak{D}_\s(\A)=\A\cdot\De,$$
where $\De=\frac{\id-\s}{g}$ and $g=\GCD\big((\id-\s)(\A)\big)$.
Furthermore, let $y\in\A$ and set
$$x=\frac{\id-\s}{g}(y)=\frac{y-\s(y)}{g}.$$
Then we have
\begin{align}\label{eq:GenWittExUFD}
    \s(g)\s(x)=\s(gx)=\s(y)-\s^2(y)=(\id-\s)(\s(y)).
\end{align} From the
definition of $g$ we know that it divides $(\id-\s)(g)=g-\s(g)$.
Thus $g$ also divides $\s(g)$. When we divide
(\ref{eq:GenWittExUFD}) by $g$ and substitute the expression for
$x$ we obtain
$$\frac{\s(g)}{g}\s\big(\frac{\id-\s}{g}(y)\big)=
\frac{\id-\s}{g}(\s(y)),$$ or, with our notation
$\De=\frac{\id-\s}{g}$,
$$\frac{\s(g)}{g}\s\big(\De(y)\big)=\De(\s(y)).$$
This shows that (\ref{eq:GenWittCond2}) holds with
\begin{align}\label{delta_g}
    \d=\s(g)/g.
\end{align}
Since $\A$ has no zero-divisors and $\s\neq\id$, it follows that
$\Ann(\De)=0$ so the equation (\ref{eq:GenWittCond1}) is clearly
true. Hence we can use Theorem \ref{thm:GenWitt} to define an
algebra structure on $\mathfrak{D}_\s(\A)=\A\cdot\De$ which
satisfies (\ref{eq:GenWittSkew}) and (\ref{eq:GenWittJacobi}) with
$\d=\s(g)/g$. Since the choice of greatest common divisor is
ambiguous (we can choose any associated element, that is, the
greatest common divisor is only unique up to a multiple by an
invertible element) this $\delta$ can be replaced by any
$\delta'=u\cdot\delta$ where $u$ is a unit (that is, an invertible
element). To see this, note that if $g'$ is another greatest
common divisor related to $g$ by $g'=u\cdot g$, then
$$\delta'=\frac{\s(ug)}{ug}=\frac{\s(u)\s(g)}{ug}=\frac{\s(u)}{u}\delta$$ and
$\s(u)/u$ is clearly a unit since $u$ is a unit. Therefore
(\ref{eq:GenWittCond2}) becomes,
\begin{align*}
    \De'=\frac{\De}{u}=\frac{\t-\s}{gu}\quad\Longrightarrow\quad
\De'(\s(a))=\frac{\s(u)}{u}\d\s(\De'(a)).
\end{align*}
\end{example}
\begin{remark}\label{rem2}
    If we choose a multiple $\Delta'=f\cdot\Delta$ of the generator
    $\Delta=\frac{\id-\s}{g}$ of $\mfD_\s(\A)$, it will generate a
    proper $\A$-submodule $\A\cdot \Delta'$ of $\mfD_\s(\A)$, unless
    $f$ is a unit. To see this, suppose on the contrary that
    $\A\cdot\Delta'=\mfD_\s(\A)$. Then there is some $g\in\A$ such
    that $g\cdot\Delta'=\Delta$. Since $\s\neq\id$
    there is some $x\in\A$ such that $\s(x)\neq x$. Then
    $$\Delta(x)=g\cdot\Delta'(x)=gf\cdot\Delta(x).$$
    Since $\Delta(x)\neq 0$ and $\A$ has no zero-divisors, we must have $gf=1$.
\end{remark}
\subsection{hom-Lie algebras} \label{subsec:homLiealg}
Let us now make the following definition.
\begin{dfn} \label{def:homLiealg} A \emph{hom-Lie algebra} $(L,\vs)$ is a
non-associative algebra $L$ together with an algebra homomorphism
$\vs :L\to L$, such that
\begin{align*}
    [x,y]_L=-[y,x]_L,
\end{align*}
\begin{align*}
    \big [(\id+\vs)(x),[y,z]_L\big ]_L+\big [(\id+\vs)(y),[z,x]_L\big ]_L+\big
[(\id+\vs)(z),[x,y]_L\big ]_L=0,
\end{align*} for all $x,y,z\in L$,
where $[\cdot,\cdot]_L$ denotes the product in $L$.
\end{dfn}
If $L$ is a Lie algebra, it is a hom-Lie algebra with its
homomorphism $\vs=\id_L$ equal to the identity operator on $L$.
\begin{example} Letting $\afrak$ be any vector space (finite- or
infinite-dimensional) we put
$$[x,y]_{\afrak}=0$$ for any $x,y\in\afrak$. Then $(\afrak, \vs_{\afrak})$ is
obviously a hom-Lie algebra for any linear map $\vs_{\afrak}$
since the above conditions are trivially satisfied. As in the Lie
case, we call these algebras \emph{abelian} or \emph{commutative}
hom-Lie algebras.
\end{example}
\begin{example}
Suppose $\A$ is a commutative associative algebra, $\s:\A\to\A$ a
homomorphism, $\De\in\mathfrak{D}_\s(\A)$ and $\d\in\A$ satisfy
the equations (\ref{eq:GenWittCond1})-(\ref{eq:GenWittCond2}).
Then since $\s(\Ann(\De))\subseteq \Ann(\De)$, the map $\s$
induces a map
$$\bar\s:\A\cdot\De\to \A\cdot\De,$$
$$\bar\s:a\cdot\De\mapsto \s(a)\cdot\De.$$
This map has the following property
\begin{align*}
    [\bar\s(a\cdot\De), \bar\s(b\cdot\De)]_\s&=[\s(a)\cdot\De, \s(b)\cdot\De]_\s=\\
    &=\big(\s^2(a)\De(\s(b))-\s^2(b)\De(\s(a))\big)\cdot\De=\\
    &=\big(\s^2(a)\d \s(\De(b))-\s^2(b)\d \s(\De(a))\big)\cdot\De=\\
    &=\d\s\big(\s(a)\De(b)-\s(b)\De(a)\big)\cdot\De=\\
    &=\d\cdot\bar\s([a\cdot\De,b\cdot\De]_\s).
\end{align*}
We suppose now that $\d\in\C\backslash\{0\}$. Dividing both sides
of the above calculation by $\d^2$ and using bilinearity of the
product, we see that
$$\frac{1}{\d}\bar\s$$
is an algebra homomorphism $\A\cdot\De\to\A\cdot\De$, and Theorem
\ref{thm:GenWitt} makes $\A\cdot\De$ with the product
$[\cdot,\cdot]_\s$ into a hom-Lie algebra with $(1/\d)\bar\s$ as
its homomorphism $\vs$.
\end{example}

By a \textit{homomorphism of hom-Lie algebras}
$\varphi:(L_1,\vs_1)\to (L_2,\vs_2)$ we mean an algebra
homomorphism from $L_1$ to $L_2$ such that $\varphi\circ
\vs_1=\vs_2\circ\varphi$, or, in other words such that the diagram
$$\begin{CD}
    L_1 @>{\varphi} >> L_2\\
    @VV{\vs_1} V  @VV{\vs_2} V\\
    L_1 @>{\varphi} >> L_2
  \end{CD}
$$
commutes. We now have the following proposition.

\begin{prop} Let $(L,\vs)$ be a hom-Lie algebra, and let $N$ be any
non-associative algebra. Let
$$\varphi:L\to N$$
be an algebra homomorphism. Then the following two conditions are
equivalent:
\begin{enumerate}
\item[\emph{1)}] There exists a linear subspace $U\subseteq N$
containing $\varphi(L)$ and a linear map $$k:U\to
    N$$ such that
    \begin{align}\label{eq:fihkfi}
        \varphi\circ\vs =k\circ\varphi.
    \end{align}
\item[\emph{2)}]
    $\ker \varphi\subseteq\ker(\varphi\circ \vs).$
\end{enumerate}
Moreover, if these conditions are satisfied, then
\begin{enumerate}
    \item[\emph{i)}] $k$ is uniquely determined on $\varphi(L)$ by $\varphi$ and $\vs$,
    \item[\emph{ii)}] $k\big|_{\varphi(L)}$ is a homomorphism
    \item[\emph{iii)}] $(\varphi(L),k\big|_{\varphi(L)})$ is a hom-Lie algebra, and
    \item[\emph{iv)}] $\varphi$ is a homomorphism of hom-Lie algebras.
\end{enumerate}
\end{prop}
\begin{remark}
    It is easy to check that condition 2) can equivalently be written
$$\vs(\ker\varphi)\subseteq\ker\varphi.$$
\end{remark}
\begin{proof}
Assume that condition 1) holds, and let $x\in\ker\varphi$. Then
$$\varphi(\vs(x))=k(\varphi(x))=k(0)=0,$$
so that $x\in\ker(\varphi\circ \vs)$. Thus 2) holds. Conversely,
assume 2) is true. Take $U=\varphi(L)$ and define $k:\varphi(L)\to
N$ by $k(\varphi(x))=\varphi(\vs(x))$. This is well-defined, since
if $\varphi(x)=\varphi(y)$ we have
$$x-y\in\ker\varphi\subseteq\ker(\varphi\circ \vs)$$ by assumption.
Hence $\varphi(\vs(x))=\varphi(\vs(y))$ so $k$ is well-defined.
Equation (\ref{eq:fihkfi}) holds by definition of $k$.

Assume now that the conditions 1) and 2) hold. To prove i), assume
that we have two linear maps $k_1:U_1\to N$ and $k_2:U_2\to N$
where $U_i$ are subspaces of $N$ with $\varphi(L)\subseteq U_i$.
Suppose they both satisfy (\ref{eq:fihkfi}). Then
$$(k_1-k_2)(\varphi(x))=\varphi(\vs(x))-\varphi(\vs(x))=0$$
for any $x\in L$. This shows that $k_1$ and $k_2$ coincide on
$\varphi(L)$. For ii) we use again the identity (\ref{eq:fihkfi}),
and that $\varphi$ is a homomorphism: (we denote the product in
$N$ by $\{\cdot,\cdot\}$ to indicate its non-associativity)
\begin{multline*}
    k(\{\varphi(x),\varphi(y)\})=k(\varphi([x,y]_L))=\varphi(\vs([x,y]_L))=
    \varphi([\vs(x),\vs(y)]_L)=\\=\{\varphi(\vs(x)),\varphi(\vs(y))\}=\{k(\varphi(x)),
    k(\varphi(y))\},
\end{multline*} for $x,y\in L$.

Using (\ref{eq:fihkfi}) and that $(L,\vs)$ is a hom-Lie algebra we
get $$\{\varphi(x),\varphi(y)\}=\varphi([x,y]_L)=\varphi(-[y,x]_L)
=-\{\varphi(y),\varphi(x)\}$$ for $x,y\in L$ and
\begin{align*}
    \circlearrowleft_{x,y,z}\{(\id+k)(\varphi(x)),\{\varphi(y),\varphi(z)\}\}&=
    \circlearrowleft_{x,y,z}\{\varphi(x)+k(\varphi(x)),\varphi([y,z]_L)\}=\\
    &=\circlearrowleft_{x,y,z}\{\varphi(x)+\varphi(\vs(x)),\varphi([y,z]_L)\}=\\
    &=\varphi(\circlearrowleft_{x,y,z}[x+\vs(x),[y,z]_L]_L)=0
\end{align*}
for $x,y,z\in L$. This shows iii), and then iv) is true since
$\varphi$ is a homomorphism satisfying (\ref{eq:fihkfi}).
\end{proof}

\subsection{Extensions of hom-Lie algebras}\label{sec:extofhlas} In this section
we will concentrate our efforts on developing the general theory
of \emph{central} extensions for hom-Lie algebras, and providing
general homological type conditions for existence of central
extensions useful for their construction.

If $U$ and $V$ are vector spaces, let $\mathrm{Alt}^2(U,V)$ denote
the space of skew-symmetric forms (alternating mappings)
$$U\times U \longrightarrow V.$$
 Exactly as in the Lie algebra case we define an extension of hom-Lie algebras with
the aid of exact sequences.
 More to the point,
 \begin{dfn} An extension of a hom-Lie algebra $(L,\vs)$ by an abelian hom-Lie
algebra $(\afrak,\vs_\afrak)$ is a commutative
 diagram with exact rows

\begin{align}\label{eq:fettcd}
    \begin{CD}
        0 @>>> \afrak @>\iota>> \hat L @>\pr>> L @>>> 0 \\
        @. @V\vs_\afrak VV @V\hat\vs VV @V\vs VV \\
        0 @>>> \afrak @>\iota>> \hat L @>\pr>> L @>>> 0 \\
    \end{CD}
\end{align}
where $(\hat L,\hat\vs)$ is a hom-Lie algebra. We say that the
extension is \emph{central} if
 $$\iota(\afrak)\subseteq Z(\hat L)=\{x\in\hat L:[x,\hat L]_{\hat L}=0\}.$$
 \end{dfn}
 The question now arises: what are the conditions for being able to construct a
central extension $\hat L$ of
 $L$? We will now derive a necessary condition for this.
The sequence above splits (as vector spaces) just as in the Lie
algebra case, meaning that there is a (linear) \emph{section}
$s:L\longrightarrow \hat L$, i.e. a linear map such that $\pr\circ
s=\id_L$. To construct a hom-Lie algebra extension we must do two
things
\begin{itemize}
   \item Define the hom-Lie algebra homomorphism $\hvs$, and
   \item construct the bracket $[\,\cdot,\,\cdot]_{\hat L}$ with the desired
properties.
\end{itemize}
Note first of all that
$$\pr\circ \hvs(x)=\vs\circ\pr(x)\quad\text{for }x\in\hat L$$ since $\pr$ is a
hom-Lie algebra homomorphism. This means that $$\pr\big
(\hvs(x)-s\circ\vs\circ\pr(x)\big )=0 $$ and this leads to, by the
exactness,
\begin{align}\label{eq:sigmahattfs}
    \hvs(x)=s\circ\vs\circ\pr(x)+\iota\circ f_s(x)
\end{align}
where $f_s: \hat L \longrightarrow\afrak$ is a function dependent on
$s$. Note that combining (\ref{eq:sigmahattfs}) with the
commutativity of the left square in (\ref{eq:fettcd}) we get for
$a\in\afrak$ that
$$\iota\circ\vs_\afrak(a)=\hat\vs\circ\iota(a)=
s\circ\vs\circ\pr\circ\iota(a)+\iota\circ
f_s\circ\iota(a)=\iota\circ f_s\circ\iota(a)$$ and hence since
$\iota$ is injective,
\begin{align}\label{eq:sigmaafrakfs}
    \vs_\afrak(a)=f_s\circ\iota(a).
\end{align}
Also $$\pr\big ([s(x),s(y)]_{\hat L}-s[x,y]_L\big )=0,$$ hence
\begin{align}\label{lift_brack}
    [s(x),s(y)]_{\hat L}=s[x,y]_L+\iota\circ g_s(x,y)
\end{align} for some $g_s\in \mathrm{Alt}^2(L,\afrak)$, a ''2-cocycle''.
This means that we have a ''lift'' of the bracket in $L$ to the
bracket in $\hat L$ for elements $x,y$ in $L$ defined by the
''2-cocycle'' and the section $s$.

Using (\ref{eq:sigmahattfs}), (\ref{lift_brack}) and the linearity
of the product we get (we temporarily suppress the indices $L$ and
$\hat L$ in the brackets and the $s$ in $g_s$), for $a,b,c\in L$,
\begin{align*}
    &\big [(\id+\hvs)(s(a)),[s(b),s(c)]\big ]=\big
[(\id+\hvs)(s(a)),s[b,c]+\iota\circ g(b,c)\big ]=\\
    &=\big [s(a)+\hvs(s(a)),s[b,c]+\iota\circ g(b,c)\big ]=\\
    &=\big [s(a),s[b,c]\big ]+\big [s(a),\iota\circ g(b,c)\big ]+\big
[\hvs(s(a)),s[b,c]\big ]+\big [\hvs(s(a)),\iota\circ g(b,c)\big
    ]=\\
    &=s\big [a,[b,c]\big ]+\iota\circ g(a,[b,c])+ \big [s(a),\iota\circ g(b,c)\big ]+\\
    &\quad+\big[s\circ\vs\circ\pr(s(a))+f_s(s(a)),s[b,c]\big]+\\
    &\quad+\big[s\circ\vs\circ\pr(s(a))+f_s(s(a)),\iota\circ g(b,c)\big]=\\
    &=s\big [a,[b,c]\big ]+\iota\circ g(a,[b,c])+ \big [s(a),\iota\circ g(b,c)\big
]+\big[s\circ\vs(a)+f_s(s(a)),s[b,c]\big]+\\
    &\quad+    \big[s\circ\vs(a)+f_s(s(a)),\iota\circ g(b,c)\big]=\\
    &=s\big [a,[b,c]\big ]+\iota\circ g(a,[b,c])+ \big [s(a),\iota\circ g(b,c)\big
]+\big[s\circ\vs(a),s[b,c]\big]
    +\\
    &\quad+\big [ f_s(s(a)),s[b,c]\big]+\big[s\circ\vs(a),\iota\circ
g(b,c)\big]+\big [f_s(s(a)),\iota\circ g(b,c)\big]=\\
    &=s\big [a,[b,c]\big ]+\iota\circ g(a,[b,c])+ \big [s(a),\iota\circ g(b,c)\big ]+\\
    &\quad+s\big[\vs(a),[b,c]\big]+\iota\circ g(\vs(a),[b,c])+\\
    &\quad+\big [ f_s(s(a)),s[b,c]\big]+ \big[s\circ\vs(a),\iota\circ
g(b,c)\big]+\big [f_s(s(a)),\iota\circ
    g(b,c)\big]=\\
    &=s\big [(\id+\vs)(a),[b,c]\big ]+\iota\circ g\big ((\id+\vs)(a),[b,c]\big )
\end{align*}where, in last step, we have used that the extension is central.
Summing up cyclically we get
\begin{align}\label{eq:twococycle}
   \cs_{a,b,c}\,g_s\big ((\id+\vs)(a),[b,c]_L\big )=0
\end{align}since $(L,\vs)$ and $(\hat{L},\hvs)$ are hom-Lie algebras.

Picking another section $\tilde{s}$, we have
$\tilde{s}(x)-s(x)=(\tilde{s}-s)(x)\in\ker\pr=\iota(\afrak)$.
Since the extension is central,
\begin{align}
    0&=[\tilde{s}(x),\tilde{s}(y)]_{\hat L}-[s(x),s(y)]_{\hat L}=\nonumber\\
    &=\tilde{s}[x,y]_L+\iota\circ g_{\tilde{s}}(x,y)-s[x,y]_L-\iota\circ
g_s(x,y)=\nonumber\\
    &=(\tilde{s}-s)([x,y]_L)+\iota\circ g_{\tilde{s}}(x,y)-\iota\circ
g_s(x,y).\label{eq:stildeg}
\end{align}
This shows that the condition (\ref{eq:twococycle}) is independent
of the section $s$. We have almost proved the following theorem.
\begin{thm}\label{thm:homext1}
    Suppose $(L,\vs)$ and $(\afrak,\vs_\afrak)$ are hom-Lie
    algebras with $\afrak$ abelian. If there exists a central extension $(\hat
L,\hat\vs)$ of $(L,\vs)$ by
    $(\afrak,\vs_\afrak)$ then for every section $s:L\to\hat L$ there is a
$g_s\in\mathrm{Alt}^2(L,\afrak)$ and a
    linear map $f_s:\hat L\to\afrak$ such that
    \begin{align}\label{eq:fandhom_a}
        f_s\circ\iota=\vs_\afrak,
    \end{align}
    \begin{align}\label{eq:bothgsandfs}
        g_s(\vs(x),\vs(y))=f_s([s(x),s(y)]_{\hat L})
    \end{align}
and
\begin{align}\label{eq:twococycle2}
    \cs_{x,y,z}\,g_s\big ((\id+\vs)(x),[y,z]_L\big )=0
\end{align} for all $x,y,z\in L$. Moreover, equation (\ref{eq:twococycle2}) is
independent of the choice of section $s$.
\end{thm}
\begin{proof}
It only remains to verify equation (\ref{eq:bothgsandfs}). We use
that $\hat\vs$ is a homomorphism. On the one hand, using
(\ref{eq:sigmahattfs}) and (\ref{lift_brack}) we have for $x,y\in
L$ that
\begin{align*}
    \hat\vs\big([s(x),s(y)]_{\hat
    L}\big)&=\hat\vs\big(s([x,y]_L)+\iota\circ g_s(x,y)\big)=\\
    &=s\circ\vs\circ\pr\circ s([x,y]_L)+ \iota\circ f_s\circ s([x,y]_L)+\\
    &\quad + s\circ\vs\circ\pr\circ\iota\circ g_s(x,y)+\iota\circ f_s\circ\iota\circ
g_s(x,y)=\\
    &=s\circ\vs\big([x,y]_L\big)+\iota\circ f_s\big([s(x),s(y)]_{\hat
    L}\big).
\end{align*}
On the other hand,
\begin{align*}
    &[\hat\vs\circ s(x),\hat\vs\circ s(y)]_{\hat L} =\\
    &=[s\circ\vs\circ\pr\circ s(x)+\iota\circ f_s\circ s(x),
    s\circ\vs\circ\pr\circ s(y)+\iota\circ f_s\circ s(y)]_{\hat L}=\\
    &=[s\circ\vs(x),s\circ\vs(y)]_{\hat L}=\\
    &=s\big([\vs(x),\vs(y)]_L\big)+\iota\circ g_s(\vs(x),\vs(y))=\\
    &=s\circ\vs\big([x,y]_L\big)+\iota\circ g_s(\vs(x),\vs(y)).
\end{align*}
Since $\iota$ is injective, (\ref{eq:bothgsandfs}) follows.
\end{proof}We now make the following definition:
\begin{dfn} A central hom-Lie algebra extension $(\hat{L},\hvs)$ of $(L,\vs)$ by
$(\afrak,\vs_{\afrak})$ is called \emph{trivial} if there exists a
linear section $s: L\longrightarrow \hat{L}$ such that
$$g_s(x,y)=0$$ for all $x,y\in L$.
\end{dfn}
\begin{remark}\label{rem:trivial}
    Note that by using (\ref{eq:stildeg}) one can show that the above definition is
equivalent to
    the statement: ''A central extension of hom-Lie algebras is trivial if and only
if for any section
    $s:L\longrightarrow\hat{L}$ there is a linear map $s_1: L\longrightarrow\hat{L}$
such that $(s+s_1)$ is a
    section and $$\iota\circ g_s(x,y)=s_1([x,y]_L)$$ for all $x,y\in L$.'' Indeed,
take a section $s$. Since the
    extension is trivial there is a section $\tilde{s}$ such that
$g_{\tilde{s}}(x,y)=0$ for all $x,y\in L$. Inserting this
    into (\ref{eq:stildeg}) gives (using that $\iota$ is one-to-one)
    $$\iota\circ g_s(x,y)=\iota\circ  g_{\tilde{s}}(x,y)+(\tilde{s}-s)[x,y]_L=
    (\tilde{s}-s)[x,y]_L$$ and putting $s_1=\tilde{s}-s$ gives necessity.
    On the other hand taking $\tilde{s}=s+s_1$ in (\ref{eq:stildeg}) gives us
sufficiency.
\end{remark}
\begin{thm}\label{thm:homext2} Suppose $(L,\vs)$ and $(\afrak,\vs_\afrak)$ are hom-Lie
algebras with $\afrak$ abelian. Then for every
$g\in\mathrm{Alt}^2(L,\afrak)$ and every linear map
$f:L\oplus\afrak\to\afrak$ such that
\begin{align}\label{eq:cfvsa}
    f(0,a)=\vs_\afrak(a) \quad\text{for } a\in\afrak,
\end{align}
\begin{align}\label{eq:cbothgsandfs}
    g(\vs(x),\vs(y))=f([x,y]_L,g(x,y))
\end{align}
and
\begin{align}\label{eq:ctwococycle2}
    \cs_{x,y,z}\,g\big ((\id+\vs)(x),[y,z]_L\big )=0,
\end{align} for $x,y,z\in L$, there exists a hom-Lie algebra $(\hat L,\hat\vs)$
which is a central extension of $(L,\vs)$ by
$(\afrak,\vs_\afrak)$.
\end{thm}
\begin{proof}
As a vector space we set $\hat L=L\oplus\afrak$. Define the
product $[\cdot,\cdot]_{\hat L}$ in $\hat L$ by setting
\begin{align}\label{eq:cproduct}
    [(x,a),(y,b)]_{\hat L}=\big([x,y]_L,g(x,y)\big)\qquad\text{for
    }(x,a),(y,b)\in\hat L
\end{align} and define $\hat\vs:\hat L\to \hat L$ by
$$\hat\vs(x,a)=(\vs(x),f(x,a))\qquad\text{for }(x,a)\in\hat L.$$
We claim that the linear map $\hat\vs$ is a homomorphism. Indeed,
\begin{align*}
    \hat\vs\big([(x,a),(y,b)]_{\hat L}\big) &= \hat\vs\big([x,y]_L,g(x,y)\big)=\\
    &=\big(\vs([x,y]_L),f([x,y]_L,g(x,y))\big)
\end{align*}
and
\begin{align*}
    [\hat\vs(x,a),\hat\vs(y,b)]_{\hat L}
    &=\big[(\vs(x),f(x,a)),(\vs(y),f(y,b))\big]=\\
    &=\big([\vs(x),\vs(y)],g(\vs(x),\vs(y))\big).
\end{align*}
These two expressions are equal because $\vs$ is a homomorphism
and (\ref{eq:cbothgsandfs}) holds. Next we prove that $(\hat
L,\hat\vs)$ is a hom-Lie algebra. Skew-symmetry of
$[\cdot,\cdot]_{\hat L}$ is immediate since $g$ is alternating.
The generalized Jacobi identity can be verified as follows:
\begin{align*}
    &\cs_{(x,a),(y,b),(z,c)}\big[(\id+\hat\vs)(x,a),\big[(y,b),(z,c)\big]_{\hat
    L}\big]_{\hat L}=\\
    &=\cs_{(x,a),(y,b),(z,c)}\big[(x+\vs(x),a+f(x,a)),([y,z]_L,g(y,z))\big]_{\hat L}=\\
    &=\cs_{(x,a),(y,b),(z,c)}\big([x+\vs(x),[y,z]_L]_L,g(x+\vs(x),[y,z]_L)\big)=0,
\end{align*}
where we used (\ref{eq:ctwococycle2}) and that $(L,\vs)$ is a
hom-Lie algebra.

Next we define $\pr$ and $\iota$ to be the natural projection and
inclusion respectively:
$$\pr:\hat L\to L,\qquad\pr(x,a)=x;$$
$$\iota:\afrak\to\hat L,\qquad \iota(a)=(0,a).$$
That the diagram (\ref{eq:fettcd}) has exact rows is now obvious.
Next we show that the linear maps $\pr$ and $\iota$ are
homomorphisms.
$$\pr\big([(x,a),(y,b)]_{\hat
L}\big)=\pr\big([x,y]_L,g(x,y)\big)=[x,y]_L=[\pr(x,a),\pr(y,b)]_L,$$
$$[\iota(a),\iota(b)]_{\hat L}=[(0,a),(0,b)]_{\hat
L}=(0,0)=\iota(0)=\iota([a,b]_\afrak)$$ since $\afrak$ was
abelian. This shows that $\pr$ and $\iota$ are homomorphisms. In
fact they are also hom-Lie algebra homomorphisms, because
$$\pr\circ\hat\vs(x,a)=\pr(\vs(x),f(x,a))=\vs(x)=\vs\circ\pr(x,a)$$
and
$$\hat\vs\circ\iota(a)=\hat\vs(0,a)=(\vs(0),f(0,a))=(0,\vs_\afrak(a))=\iota\circ\vs_\afrak(a),$$
where we used (\ref{eq:cfvsa}). This proves that $(\hat
L,\hat\vs)$ is an extension of $(L,\vs)$ by $(\afrak,\vs_\afrak)$.
Finally, that the extension is central is clear from the
definition of $\iota$ and (\ref{eq:cproduct}).
\end{proof}

\section{Examples}
\subsection{A $q$-deformed Witt algebra} \label{sec:qwitt}
Let $\A$ be the complex algebra of Laurent polynomials in one
variable $t$, i.e.
$$\A=\C[t,t^{-1}]\cong \,
\C[x,y]\big /(xy-1).$$ Fix $q\in \C\setminus\{0,1\}$, and let $\s$
be the unique endomorphism on $\A$ determined by
$$\s(t)=qt.$$
Explicitly, we have
$$\s(f(t))=f(qt),\qquad\text{for } f(t)\in\A.$$
The set $\mfD_{\s}(\A)$ of all $\s$-derivations on $\A$ is a free
$\A$-module of rank one, and the mapping
$$D:\A\to\A, $$ defined by
\begin{align} \label{eq:jacksonq}
    D(f(t))= t\frac{\s(f(t))-f(t)}{\s(t)-t}=\frac{f(qt)-f(t)}{q-1}\qquad\text{for
}f(t)\in\A,
\end{align} is a generator.  \\

To see that $D$ indeed generates $\mfD_{\s}(\A)$, note that, since
$\A$ is a UFD, a generator of $\mfD_{\s}(\A)$ is on the form
$$\frac{\id-\s}{\GCD\big ((\id-\s)(\A)\big )},$$ by Theorem
\ref{thm:polndifferent}. Now, a greatest common divisor on
$(\id-\s)(\A)$ is any element of $\A$ on the form $ct^k$, where
$c\in \C\setminus\{0\}$ and $k\in\Z$. This is because
$\gcd\big((\id-\s)(\A)\big)$ divides any element of $(\id-\s)(\A)$
so in particular it divides $(\id-\s)(t)=-(q-1)t$ which is a unit
(when $q\neq 1$). This means that $q-1$ is a $\GCD\big
((\id-\s)(\A)\big )$. Therefore,
$$D(f(t))=\frac{f(qt)-f(t)}{q-1}=-\frac{\id-\s}{q-1}\big (f(t)\big )$$ and so
$D=-\frac{\id-\s}{q-1}$ is a generator for $\mfD_\s(\A)$.
\begin{remark} Note that $t^{-1}D=D_q$, the \textit{Jackson $q$-derivative}.
\end{remark}

 Since $D$ is a polynomial (over $\C$) in $\s$, $D$ and $\s$ commute.
Let $[\cdot,\cdot]_\s$ denote the product on $\mfD_{\s}(\A)$
defined by
\begin{align}\label{eq:bracketdef}
    [f\cdot D, g\cdot D]_\s=(\s(f)\cdot D)\circ(g\cdot D)- (\s(g)\cdot
D)\circ(f\cdot D)
\end{align} for $f,g\in\A$ coming from (\ref{eq:GenWittProdDef}). It satisfies the
following identities:
\begin{align} \label{eq:bracketidentity}
[f\cdot D, g\cdot D]_\s =\big(\s(f)D(g)-\s(g)D(f)\big)\cdot D,
\end{align}
\begin{align} \label{eq:anticommutativity}
    [f\cdot D, g\cdot D]_\s = -[g\cdot D, f\cdot D]_\s,
\end{align}
and
\begin{align} \label{eq:jacobi}
    \big [(\s(f)+f)\cdot D, [g\cdot D, h\cdot D]_\s\big ]_\s&+\big [(\s(g)+g)\cdot
D, [h\cdot D, f\cdot D]_\s\big ]_\s
    +\notag\\
    &+\big [(\s(h)+h)\cdot D, [f\cdot D, g\cdot D]_\s\big ]_\s=0,
\end{align}
for all $f,g,h\in\A$. The identities (\ref{eq:anticommutativity})
and (\ref{eq:jacobi}) show that $\mfD_{\s}(\A)$ is a hom-Lie
algebra with
$$\vs:\mfD_{\s}(\A)\to\mfD_{\s}(\A)$$
$$\vs:f\cdot D\mapsto \s(f)\cdot D$$
as its homomorphism. As a $\C$-linear space, $\mfD_{\s}(\A)$ has a
basis $\{d_n\,|\,n\in\Z\}$, where
\begin{align}\label{eq:qwittbasis}
    d_n=-t^n\cdot D.
\end{align} Note that
\begin{align} \label{eq:sigmaOnBasisElement}
    \s(-t^n)=-q^nt^n,
\end{align} which imply
\begin{align} \label{eq:hOnBasisElement}
    \vs(d_n)=q^nd_n.
\end{align} Note further that
\begin{align} \label{eq:jacksonOnBasisElement}
    D(-t^n)=\frac{-q^nt^n+t^n}{q-1}= -\{n\}t^n,
\end{align} where $\{n\}$ for
$n\in\Z$ denotes the \emph{$q$-number}
$$\{n\}=\frac{q^n-1}{q-1}.$$ Using (\ref{eq:bracketidentity}) with $f(t)=-t^n$ and
$g(t)=-t^l$ we obtain the following important commutation
relation:
\begin{align} \label{eq:qwittcomrel}
    [d_n,d_l]_\s&=\big((-q^nt^n)\cdot(-\{l\}t^l))-
    (-q^lt^l)\cdot(-\{n\}t^n)\big)\cdot D=\notag\\
    &=\big(q^l\frac{q^n-1}{q-1}-q^n\frac{q^l-1}{q-1}\big)\cdot
    (-t^{n+l})\cdot D)=\notag\\
    &=\frac{q^n-q^l}{q-1}d_{n+l}=\big(\{n\}-\{l\}\big)d_{n+l},
\end{align}
for $n,l\in\Z$, where the bracket is defined on generators by
(\ref{eq:GenWittProdDef}) as
$$[d_n,d_l]_\s=q^n d_nd_l-q^l d_ld_n.$$ This means, in particular, that
$\mfD_\s(\A)$ admits a $\Z$-grading as an algebra:
$$\mfD_\s(\A)=\bigoplus_{i\in\Z}\C\cdot d_i.$$
\begin{remark} If $q=1$ we simply define $D$
to be $t\cdot\partial$ where $\partial=\frac{d}{dt}$, the ordinary
differential operator. Then $\partial$ generates $\mfD_{\id}(\A)$
even though Theorem \ref{thm:polndifferent} cannot be used. The
relation
$$[d_n,d_l]_\s=\big(\{n\}-\{l\}\big)d_{n+l}$$ then becomes the standard commutation
relation in the Witt algebra:
$$[\partial_n,\partial_l]=(n-l)\partial_{n+l},$$
where $\partial_n=-t^n\cdot D$.
\end{remark}
It follows from (\ref{eq:anticommutativity}) that
\begin{align} \label{eq:qwittskew}
    [d_n,d_l]_\s=-[d_l,d_n]_\s,
\end{align} and substituting $f(t)=-t^n$, $g(t)=-t^l$ and $h(t)=-t^m$ into
(\ref{eq:jacobi}) we obtain the following $q$-deformation of the
Jacobi identity:
\begin{align} \label{eq:qwittjacobi}
    (q^n+1)\big [d_n, [d_l, d_m]_\s\big ]_\s&+ (q^l+1)\big [d_l,  [d_m, d_n
]_\s\big]_\s+\notag\\
    &+(q^m+1)\big [d_m, [d_n,d_l]_\s\big ]_\s\quad=0,
\end{align} for all $n,l,m\in\Z$. Hence,
\begin{thm} \label{th:qWittlinear} Let $\A=\C[t,t^{-1}]$, Then the $\C$-linear space
    $$\mfD_\s (\A) = \bigoplus_{n\in\Z} \C\cdot d_n,$$
    where
    $$D=-\frac{\id-\s}{q-1},$$ $d_n=-t^nD$ and $\s(t)=qt$ can be equipped with the
    bracket multiplication
    $$[\,\cdot,\,\cdot]_\s\,:\, \mfD_\s (\A)\times\mfD_\s (\A) \longrightarrow
    \mfD_\s (\A)$$ defined on generators by (\ref{eq:GenWittProdDef}) as
    $$[d_n,d_m]_\s=q^nd_{n}d_m-q^md_{m}d_n$$ with commutation relations
    $$[d_n,d_m]_\s=\big(\{n\}-\{m\}\big)d_{n+m}.$$
    This bracket satisfies skew-symmetry
    $$[d_n,d_m]_\s=-[d_m,d_n]_\s,$$ and a $\s$-deformed Jacobi-identity
    \begin{align*}
        (q^n+1)\big [d_n, [d_l, d_m]_\s\big ]_\s&+ (q^l+1)\big [d_l,  [d_m, d_n
]_\s\big]_\s+\notag\\
    &+(q^m+1)\big [d_m, [d_n,d_l]_\s\big ]_\s\quad=0.
    \end{align*}
\end{thm}

\begin{remark}
The associative algebra with an infinite number of (abstract)
generators $\{d_j \mid j \in \mathbb{Z} \}$ and defining relations
$$q^nd_{n}d_m-q^md_{m}d_n=\big(\{n\}-\{m\}\big)d_{n+m}, \quad n,m \in \Z$$
is a well defined associative algebra, since our construction,
summarized in Theorem \ref{th:qWittlinear}, yields at the same
time its operator representation. Naturally, an outcome of our
approach is that this parametric family of algebras is a
deformation of the Witt algebra defined by relations
\eqref{rel:Witt} in the sense that \eqref{rel:Witt} is obtained
when $q=1$.
\end{remark}

\subsection{Non-linearly deformed Witt algebras}
With the aid of Theorems \ref{thm:polndifferent} and \ref{thm:GenWitt} we will now construct a non-linear
deformation of the derivations of $\A=\C[t,t^{-1}]$, the algebra of Laurent polynomials. Take any $p(t)\in\A$
and assume that $\s(t)=p(t)$. In addition, we assume $\s(1)=1$, since if this is not the case, we would have had
$\s(1)=0$ because $\A$ has no zero-divisors and so $\sigma(1)=0$ would imply $\s=0$ identically. This leads us
to
$$1=\s(1)=\s(t\cdot t^{-1})=\s(t)\s(t^{-1})\Longrightarrow \s(t^{-1})=\s(t)^{-1},$$
implying two things:
\begin{enumerate}
    \item $\s(t)$ must be a unit, and
    \item $\s(t^{-1})$ is completely determined by $\s(t)$ as its inverse in $\A$.
\end{enumerate}Hence, since $\s(t)$ is a unit, $\s(t)=p(t)=qt^s$, for some
$q\in\C\setminus\{0\}$ and $s\in\Z$. We will, however, continue
writing $p(t)$ instead of $qt^s$ except in the explicit
calculations.\\

It suffices to compute a greatest common divisor of $(\id-\s)(\A)$ on the generator $t$ since $\s(t^{-1})$ is
determined by $\s(t)$. Furthermore, any $\gcd$ is only determined up to a multiple of a unit. This gives us that
$$g=\alpha^{-1}t^{k-1}(\id-\s)(t)=\alpha^{-1}t^{k-1}(t-p(t))=\alpha^{-1} t^{k-1}(t-qt^{s})$$ is a perfectly general $\gcd$ and so Theorem
\ref{thm:polndifferent} tells us that
\begin{align*}
    D=\frac{\id-\s}{\alpha^{-1}t^{k-1}(t-p(t))}=\alpha t^{-k+1}\frac{\id-\sigma}{t-qt^{s}}=
    \alpha t^{-k}\frac{\id-\sigma}{1-qt^{s-1}}
\end{align*} is a generator for $\mfD_\s(A)$.\\

Two direct consequences of this is that, firstly, if $r\in\Z_{\geq 0}$, then
\begin{align*}
    D(t^r)&=\alpha\cdot t^{-k+r}\frac{1-q^rt^{r(s-1)}}{1-qt^{s-1}}=\alpha\sum_{l=0}^{r-1}q^lt^{(s-1)l+r-k}=\\
    &=\alpha t^{-k}\sum_{l=0}^{r-1}p(t)^lt^{r-l}=\alpha t^{-k}\sum_{l=0}^{r-1}p(t)^{r-1-l}t^{l+1}
\end{align*} and secondly, if $r\in\Z_{<0}$, then
\begin{align*}
    D(t^r)&=\alpha\cdot t^{-k+r}\frac{1-q^rt^{r(s-1)}}{1-qt^{s-1}}=-\alpha\cdot
    t^{-k+r}q^rt^{r(s-1)}\frac{1-q^{-r}t^{-r(s-1)}}{1-qt^{s-1}}=\\
    &=-\alpha\cdot
    t^{-k+r}q^rt^{r(s-1)}\sum_{l=0}^{-r-1}q^lt^{l(s-1)}=-\alpha\sum_{l=0}^{-r-1}q^{r+l}t^{(r+l)(s-1)-k+r}.
\end{align*}

The $\s$-derivations on $\C[t,t^{-1}]$ are on the form $f(t)\cdot
D$ for $f\in\C[t,t^{-1}]$ and so, given that $t^{\Z}$ is a linear
basis of $\C[t,t^{-1}]$ (over $\C$), $-t^{\Z}\cdot D$ is a linear
basis (over $\C$ again) for $\mfD_\s(\C[t,t^{-1}])$. We now
introduce a bracket on $\mfD_\s(\C[t,t^{-1}])$ in accordance with
Theorem \ref{thm:GenWitt} as we did in the previous section. Once
again,
\begin{align*}
    [-t^n\cdot D,-t^m\cdot D]_\sigma=\big (\sigma(-t^n) D(-t^m)-\sigma(-t^m) D(-t^n)\big ) D.
\end{align*} To continue we consider three cases
(1) $n,m>0$, (2) $n>0$, $m<0$, and (3) $n,m<0$.
\subsubsection*{Case 1.} Assume $n,m>0$. Thus the coefficient in the bracket is
    \begin{align*}
        \sigma(t^n)D(t^m)&-\sigma(t^m)D(t^n)=p(t)^n\cdot\alpha t^{-k}\sum^{m-1}_{l=0}p(t)^{m-1-l}t^{l+1}-\\
        &-p(t)^m\cdot\alpha t^{-k}\sum_{l=0}^{n-1}p(t)^{n-1-l}t^{l+1}=\\
        &=\alpha\Big (\sum^{m-1}_{l=0}p(t)^{n+m-1-l}t^{l-k+1}-\sum_{l=0}^{n-1}p(t)^{n+m-1-l}t^{l-k+1}\Big ).
    \end{align*}To re-write this we use the ''sign function''
    $$\sign(x)=\begin{cases}
    -1 & \text{if } x<0\\
    0 & \text{if } x=0\\
    1 & \text{if } x>0.
    \end{cases}$$So,
    \begin{align*}
        \alpha\Big (&\sum^{m-1}_{l=0}p(t)^{n+m-1-l}t^{l-k+1}-\sum_{l=0}^{n-1}p(t)^{n+m-1-l}t^{l-k+1}\Big )=\\
        &=\alpha\sign(m-n)\sum_{l=\min(n,m)}^{\max(n,m)-1}p(t)^{n+m-1-l}t^{l-k+1}=\\
        &=\alpha\sign(m-n)\sum_{l=\min(n,m)}^{\max(n,m)-1}q^{n+m-1-l}t^{(n+m-1)s-(s-1)l-(k-1)}
    \end{align*} giving that (for $n,m\in\Z_{\geq 0}$)
    \begin{align}\label{bracket:nonlinear}
    [d_n,d_m]_\sigma=\alpha\sign(n-m)\sum_{l=\min(n,m)}^{\max(n,m)-1}q^{n+m-1-l}d_{(n+m-1)s-(s-1)l-(k-1)}.
    \end{align}
    \begin{remark}\label{rem:nonlinear}
    Note that if we take $k=0,\,s=1$ and $\alpha=1$, the right hand sum in
(\ref{bracket:nonlinear}) contains only
    the generator $d_{n+m}$ multiplied by the coefficient
    \begin{align*}
        &\sign(n-m)q^{n+m-1}\sum_{l=\min(n,m)}^{\max(n,m)-1}(q^{-1})^l=\\
        &=\sign(n-m)q^{n+m-1}\frac{(q^{-1})^{\max(n,m)}-
        (q^{-1})^{\min(n,m)}}{q^{-1}-1}=\\
        &=\sign(n-m)\frac{q^{\max(n,m)}-q^{\min(n,m)}}{q-1}=\\
        &=\sign(n-m)\big (\{\max(n,m)\}-\{\min(n,m)\}\big )=\{n\}-\{m\}.
    \end{align*}So the commutation relation (\ref{bracket:nonlinear}) reduces to the
    relation
    (\ref{eq:qwittcomrel}) for the $q$-deformation of the Witt algebra described in
    section 2.1.
\end{remark}
\subsubsection*{Case 2.} Now, suppose $n>0$ and $m<0$. We then get the bracket coefficient
    \begin{align*}
        \sigma&(t^n)D(t^m)-\sigma(t^m)D(t^n)=\\
        &=-\alpha\cdot q^nt^{ns}\sum_{l_1=0}^{-m-1}q^{m+l_1}t^{(m+l_1)(s-1)-k+m}-\alpha\cdot
        q^mt^{ms}\sum_{l_2=0}^{n-1}q^{l_2}t^{(s-1)l_2+n-k}=\\
        &=-\alpha\Big (\sum_{l_1=0}^{-m-1}q^{n+m+l_1}t^{(m+l_1)(s-1)-k+m+ns}+
        \sum_{l_2=0}^{n-1}q^{m+l_2}t^{(s-1)l_2+n-k+ms}\Big ).
    \end{align*}We now show that there is no overlap between these two sums. Knowing that
    $0\leq l_2\leq n-1$ we consider the difference in exponents of $t$:
    \begin{multline*}
        (m+l_1)(s-1)-k+m+ns-(s-1)l_2-n+k-ms=\\
        =(s-1)(l_1-l_2)+m(s-1)+(1-s)m+(s-1)n=(s-1)(l_1-l_2+n)
    \end{multline*}and this is zero (for $s\neq 1$) when
    $n=l_2-l_1$. But, $n=l_2-l_1\leq n-1-0=n-1$ which is a contradiction and hence we cannot have any overlap.
    Hence, we see that the bracket becomes
    \begin{align}\label{bracket:nonlinear2}
    [d_n,d_m]_\sigma =\alpha\Big (\sum_{l=0}^{-m-1}q^{n+m+l}d_{(m+l)(s-1)+ns+m-k}+
        \sum_{l=0}^{n-1}q^{m+l}d_{(s-1)l+n+ms-k}\Big ).
    \end{align}

\subsubsection*{Case 2$^{\prime}$.} By interchanging the role of $n$ and $m$ so
that $m>0$ and $n<0$ we get instead
\begin{align*}
        \sigma&(t^n)D(t^m)-\sigma(t^m)D(t^n)=\\
       &=\alpha\Big (\sum_{l_1=0}^{m-1}q^{n+l_1}t^{(s-1)l_1+m+ns-k}+
         \sum_{l_2=0}^{-n-1}q^{m+n+l_2}t^{(n+l_2)(s-1)+n+ms-k}\Big ),
    \end{align*}so the coefficient for the bracket becomes
\begin{align}\label{bracket:nonlinear2prime}
     [d_n,d_m]_\sigma&=\alpha\Big (\sum_{l_1=0}^{m-1}q^{n+l_1}
        d_{(s-1)l_1+m+ns-k}+
         \sum_{l_2=0}^{-n-1}q^{m+n+l_2}d_{(n+l_2)(s-1)+n+ms-k}\Big ).
\end{align}
\begin{remark}If we put $k=0$, $s=1$ and $\alpha=1$ in Case 2 and $2^{\prime}$, we once again get the single generator $d_{n+m}$ multiplied
with
\begin{align*}
    q^{n+m}\sum_{l=0}^{-m-1}q^l+q^m\sum_{l=0}^{n-1}q^l&=q^{m+n}\frac{1-q^{-m}}{1-q}+q^m\frac{1-q^n}{1-q}=\\
    &=-q^n\frac{1-q^m}{1-q}+q^m\frac{1-q^n}{1-q}=\frac{q^m-q^n}{1-q}=\\
    &=\{n\}-\{m\},
\end{align*}just as we would expect from the case in the previous section.
\end{remark}
\subsubsection*{Case 3.} Both $n,m<0$. This leads to
\begin{align*}
    \sigma&(t^n)D(t^m)-\sigma(t^m)D(t^n)= \\
     &=-\alpha\cdot
q^nt^{ns}\sum_{l_1=0}^{-m-1}q^{m+l_1}t^{(m+l_1)(s-1)-k+m}+\\
&\qquad+\alpha\cdot q^mt^{ms}\sum_{l_2=0}^{-n-1}q^{n+l_2}t^{(n+l_2)(s-1)-k+n}=\\
    &=-\alpha\Big (\sum_{l_1=0}^{-m-1}q^{n+m+l_1}t^{(m+n)s+(s-1)l_1-k}-
        \sum_{l_2=0}^{-n-1}q^{n+m+l_2}t^{(n+m)s+(s-1)l_2-k}\Big ).
\end{align*} This leads to a bracket coefficient resembling that of case 1, namely
\begin{align}\label{bracket:nonlinear3}
    [d_n,d_m]_\sigma=\alpha\sign(n-m)\sum_{l=\min(-n,-m)}^{\max(-n,-m)-1}q^{n+m+l}d_{(m+n)s+(s-1)l-k}.
\end{align}

We can now, from (\ref{delta_g}), calculate $\d$ to get
\begin{align*}
    \delta&=
    \frac{\sigma(g)}{g}=\frac{\alpha^{-1}t^{sk}(1-q^st^{s(s-1)})}
{\alpha^{-1}t^k(1-qt^{s-1})}=q^kt^{k(s-1)}\sum_{r=0}^{s-1}(qt^{s-1})^r.
\end{align*}
 This means, by the
definition of $\d$, that $D$ and $\s$ span a ''quantum
plane''-like commutation relation
$$D\s=q^kt^{k(s-1)}\sum_{r=0}^{s-1}(qt^{s-1})^r \cdot \s D.$$
To get a hom-Lie algebra it is enough for $\d$ to be a (non-zero)
complex number and this can be achieved only when $s=1$, that is,
when the deformation is linear
(i.e. when $\s$ homogeneous of degree zero).\\

Theorem \ref{thm:GenWitt} now tells us what a generalized Jacobi
identity looks like
\begin{align*}
    \cs_{n,m,l}\,\Big (q^n\big [d_{ns},[d_m,d_l]_\s\big
]_\s+q^kt^{k(s-1)}\sum_{r=0}^{s-1}(qt^{s-1})^r\big [
     d_n[d_m,d_l]_\s\big ]_\s\Big )=0.
\end{align*}The hom-Lie algebra Jacobi-identity ($s=1$) becomes
\begin{align*}
    \cs_{n,m,l}\,(q^n+q^k)\big [d_{n},[d_m,d_l]_\s\big ]_\s=0.
\end{align*}
We summarize our findings in a theorem.
\begin{thm} \label{th:qWittnonlinear} Let $\A=\C[t,t^{-1}]$. Then the $\C$-linear space
    $$\mfD_\s(\A) = \bigoplus_{n\in\Z} \C\cdot d_n,$$
    where
    $$D=\alpha t^{-k+1}\frac{\id-\s}{t-qt^s},$$ $d_n=-t^nD$ and $\s(t)=qt^s$, can be
equipped
    with the bracket product
    $$[\,\cdot,\,\cdot]_\s\,:\, \mfD_\s(\A)\times \mfD_\s(\A)\longrightarrow
    \mfD_\s(\A)$$ defined on generators by (\ref{eq:GenWittProdDef}) as
    $$[d_n,d_m]_\s=q^nd_{ns}d_m-q^md_{ms}d_n$$ and satisfying
    defining commutation relations
        \begin{align*}
            &[d_n,d_m]_\sigma=\alpha\sign(n-m)\sum_{l=\min(n,m)}^{\max(n,m)-1}q^{n+m-1-l}
            d_{s(n+m-1)-(k-1)-l(s-1)}\\
            &\qquad\text{for } n,m\geq 0;\\
            &[d_n,d_m]_\sigma=\alpha\Big
            (\sum_{l=0}^{-m-1}q^{n+m+l}d_{(m+l)(s-1)+ns+m-k}+\sum_{l=0}^{n-1}q^{m+l}d_{(s-1)l+n+ms-k}\Big )\\
            &\qquad\text{for } n\geq 0, m<0;\\
            &[d_n,d_m]_\sigma=\alpha\Big (\sum_{l_1=0}^{m-1}q^{n+l_1}
        d_{(s-1)l_1+m+ns-k}+
         \sum_{l_2=0}^{-n-1}q^{m+n+l_2}d_{(n+l_2)(s-1)+n+ms-k}\Big )\\
            &\qquad\text{for } m\geq 0, n<0;\\
            &[d_n,d_m]_\sigma=\alpha\sign(n-m)\sum_{l=\min(-n,-m)}^{\max(-n,-m)-1}q^{n+m+l}d_{(m+n)s+(s-1)l-k}\\
            &\qquad\text{for } n,m<0.
        \end{align*}
    Furthermore, this bracket satisfies skew-symmetry
    $$[d_n,d_m]_\s=-[d_m,d_n]_\s,$$ and a $\s$-deformed Jacobi-identity
    \begin{align*}
        \cs_{n,m,l}\,\Big (q^n\big [d_{ns},[d_m,d_l]_\s\big
]_\s+q^kt^{k(s-1)}\sum_{r=0}^{s-1}(qt^{s-1})^r\big [
        d_n[d_m,d_l]_\s\big ]_\s\Big )=0.
    \end{align*}
\end{thm}

\begin{remark}
The associative algebra with an infinite number of (abstract)
generators $\{d_j \mid j \in \mathbb{Z} \}$ and defining relations
$$q^nd_{ns}d_m-q^md_{ms}d_n=\begin{cases}
    \text{Eq. (\ref{bracket:nonlinear})} & \text{for } n,m\geq 0\\
    \text{Eq. (\ref{bracket:nonlinear2})} & \text{for } n\geq 0, m<0\\
    \text{Eq. (\ref{bracket:nonlinear2prime})} & \text{for } m\geq 0, n<0\\
    \text{Eq. (\ref{bracket:nonlinear3})} & \text{for } n,m<0
  \end{cases}
$$
is a well defined associative algebra, since our construction yields at the same time its operator
representation. Naturally, an outcome of our approach is that this parametric family of algebras is a
deformation of the Witt algebra defined by relations \eqref{rel:Witt} in the sense that \eqref{rel:Witt} is
obtained when $q=1$ and $k=0$, $s=1$, $\alpha = 1$.
\end{remark}

\subsubsection{A submodule of $\mfD_\s(\C[t,t^{-1}])$} We let, as before,
$\A=\C[t,t^{-1}]$, the algebra of Laurent polynomials, and $\s$ be
some non-zero endomorphism such that $\s(t)=p(t)$. In the previous
section we showed that any greatest common divisor of
$(\id-\s)(\A)$ has the form
$\alpha^{-1}t^{k-1}(t-p(t))=\alpha^{-1}t^{k}(1-qt^{s-1})$ for
$k\in\Z$ and nonzero $\alpha\in\C$. As described in Remark
\ref{rem2}, this means that
$$D=\frac{\id-\s}{\alpha^{-1}\cdot t^k}$$ generates a proper cyclic $\A$-submodule
$\M$ of  $\mfD_\s(\A)$,
unless $p(t)=\beta t$ for some $\beta\in\C$.\\

As above we calculate $\d$ using (\ref{delta_g}) and we find
$$\d=q^kt^{(s-1)k}$$ which means that $D$ and $\s$ satisfy the following relation
$$D\s=q^kt^{(s-1)k}\s D.$$
We set
$$\dh_n=-t^n D.$$
Before we calculate the bracket, we note that $D(t)=\a
t^{-k+1}(1-qt^{s-1})$, and
$$\a\frac{\id-\s}{t^k}\big (t^r\big )=\a t^{r-k}(1-q^rt^{(s-1)r}).$$
The coefficient of $D$ in the bracket $[\dh_n,\dh_m]_\s$ then
becomes
\begin{align*}
    \s(t^n)Dt^m-\s(t^m)Dt^n&=\a q^nt^{ns+m-k}(1-q^mt^{(s-1)m})-\\
    &\quad-\a q^mt^{ms+n-k}(1-q^nt^{(s-1)n})-\\
    &=\a q^nt^{ns+m-k}-\a q^mt^{ms+n-k}
\end{align*} which means that
\begin{align*}
    [\dh_n,\dh_m]_\s=\a q^m\dh_{ms+n-k}-\a q^n\dh_{ns+m-k},
\end{align*}where
$$[\dh_n,\dh_m]_\s=q^n\dh_{ns}\dh_m-q^m\dh_{ms}\dh_n$$ by (\ref{eq:GenWittProdDef}).
Putting $a=-t^n$, $b=-t^m$ and $c=-t^l$ in
(\ref{eq:GenWittJacobi}) with $\delta=q^kt^{(s-1)k}$ and $\De=D$,
we get
\begin{multline*}
    \cs_{n,m,l}\,\Big (\big [-q^nt^{ns} D,[-t^m D, -t^l D]_\s\big
]_\s+q^kt^{(s-1)k}\big [-t^n D,[-t^m D,-t^l D]_\s\big
        ]_\s\Big )=\\
    =\cs_{n,m,l}\,\Big (q^n\big [\dh_{ns},[\dh_m,\dh_l]_\s\big
]_\s+q^kt^{(s-1)k}\big [\dh_n,[\dh_m,\dh_l]_\s\big
    ]_\s\Big )=0.
\end{multline*}
We summarize the obtained results in the following theorem.
\begin{thm} \label{th:submodWitt1} The $\C$-linear space
        $$\M=
        \bigoplus_{i\in\Z}\C\cdot\dh_i\quad \text{ with }\quad \dh_i=-t^i D$$ of
$\s$-derivations on
        $\A$ allows a structure as an algebra with
        bracket defined on generators (by (\ref{eq:GenWittProdDef})) as
        $$[\dh_n,\dh_m]_\s=q^n\dh_{ns}\dh_m-q^m\dh_{ms}\dh_n$$ and
        satisfying relations
        $$[\dh_n,\dh_m]_\s= \alpha q^m\dh_{ms+n-k}-\alpha q^n\dh_{ns+m-k},$$ with
$s\in\Z$ and $\alpha\in\C$.
        The $\s$-deformed Jacobi identity becomes
        \begin{align*}
        \cs_{n,m,l}\,\Big (q^n\big [\dh_{ns},[\dh_m,\dh_l]_\s\big
]_\s+q^kt^{(s-1)k}\big [\dh_n,[\dh_m,\dh_l]_\s\big
        ]_\s\Big )=0.
        \end{align*}
\end{thm}
\begin{remark}
    The only possible way to obtain a $\Z$-grading on $\M$ with the bracket
$[\cdot,\cdot]_\s$ is
    when $k=0$ and $s=1$ in the above theorem.
\end{remark}
\begin{remark}
    By performing a change of basis and consider instead $\dh_n=-t^{n+k} D$
    in the definition of $\dh_n$
    we can evade the use of $k$ altogether.
    Hence we see that the $k$-shifted grading is something resulting
    from a choice of basis for $\M$.
\end{remark}
\begin{remark}
The associative algebra with an infinite number of (abstract)
generators $\{d_j \mid j \in \mathbb{Z} \}$ and defining relations
$$q^n d_{ns} d_m-q^m d_{ms}d_n=
\alpha q^m d_{ms+n-k}-\alpha q^n d_{ns+m-k}, \quad n,m \in \Z$$ is
a well defined associative algebra, since our construction,
summarized in Theorem \ref{th:submodWitt1}, yields at the same
time its operator representation. It is interesting that, when
$q=1$, $k=0$ and $s=1$, we get a commutative algebra with
countable number of generators instead of the Witt algebra.
\end{remark}

\subsubsection{Generalization to several variables}
We let the boldface font denote an $\N$-vector, e.g.
$$\kbf=(k_1,k_2,\dots,k_n),\qquad k_i\in\N.$$ Consider the algebra of Laurent
polynomials in $z_1,z_2,\dots, z_n$
$$\A=\C[z_1^{\pm 1},z_2^{\pm 1},\dots,z_n^{\pm 1}]\cong
\frac{\C[z_1,\dots,z_n,u_1,\dots,u_n]}{(z_1u_1-1,\dots,z_nu_n-1)}$$
and let
$$\s(z_1)=q_{z_1}z_1^{S_{1,1}}\cdots
z_n^{S_{1,n}},\quad\dots\quad,\,\s(z_n)=q_{z_n}z_1^{S_{n,1}}\cdots
z_n^{S_{n,n}}.$$ Notice that $\s$ is determined by a matrix
$$S=[S_{ij}]$$ and the complex numbers $q_{z_k}$. A common divisor
on $(\id-\s)(\A)$ is $$g=Q^{-1}z_1^{G_1}\cdots z_n^{G_n}$$ and so
the element
$$D=Q\frac{\id-\s}{z_1^{G_1}\cdots z_n^{G_n}}$$ generates an $\A$-submodule $\M$ of
$\mfD_\s(\A)$. Using these $D$
 and $\s$ we calculate $\d$ to be (by formula (\ref{eq:GenWittCond2}))
\begin{align*}
    \d&=q_{z_1}^{G_1}\cdots
q_{z_n}^{G_n}z_1^{(S_{1,1}-1)G_1+S_{2,1}G_2+\cdots+S_{n,1}G_n}\cdots
        z_n^{S_{1,n}G_1+S_{2,n}G_2+\cdots+(S_{n,n}-1)G_n}=\\
        &=q_{z_1}^{G_1}\cdots q_{z_n}^{G_n}z_1^{\d_1}\cdots z_n^{\d_n},
\end{align*}where,
$$\d_k=S_{1,k}G_1+S_{2,k}G_2+\cdots+(S_{k,k}-1)G_k+\cdots +S_{n,k}G_n.$$ We also
introduce the following notation
$$\alpha_r(\lbf)=\sum_{i=1}^{n}S_{i,r}l_i.$$ Now,
\begin{align*}
    \s(z_1)^{l_1}\cdots \s(z_n)^{l_n}&=q_{z_1}^{l_1}\cdots
q_{z_n}^{l_n}z_1^{S_{1,1}l_1+\cdots+S_{n,1}l_n}
    \cdots z_n^{S_{1,n}l_1+\cdots+S_{n,n}l_n}=\\
    &=q_{z_1}^{l_1}\cdots q_{z_n}^{l_n}z_1^{\alpha_1(\lbf)}\cdots z_n^{\alpha_n(\lbf)}
\end{align*} and so,
\begin{align*}
    D(z_1^{l_1}\cdots z_n^{l_n})&=Q\frac{\id-\s}{z_1^{G_1}\cdots
z_n^{G_n}}(z_1^{l_1}\cdots z_n^{l_n})=\\
    &=Q\frac{z_1^{l_1}\cdots
z_n^{l_n}-\s(z_1)^{l_1}\cdots\s(z_n)^{l_n}}{z_1^{G_1}\cdots z_n^{G_n}}=\\
    &=Q\frac{z_1^{l_1}\cdots z_n^{l_n}-q_{z_1}^{l_1}\cdots
q_{z_n}^{l_n}z_1^{\alpha_1(\lbf)}\cdots
    z_n^{\alpha_n(\lbf)}}{z_1^{G_1}\cdots z_n^{G_n}}=\\
    &=-Qz_1^{l_1-G_1}\cdots z_n^{l_n-G_n}\Big (q_{z_1}^{l_1}\cdots
q_{z_n}^{l_n}z_1^{\alpha_1(\lbf)-l_1}\cdots
        z_n^{\alpha_n(\lbf)-l_n}-1\Big )
\end{align*}We now put
$$d_{\lbf}=d_{l_1,\dots, l_n}=-z_1^{l_1}\cdots z_n^{l_n}D$$ and calculate the
coefficient of the bracket $[d_{\kbf},d_{\lbf}]_\s$ with the aid
of Theorem \ref{thm:GenWitt} as before:
\begin{multline*}
    \s(z_1^{k_1}\cdots z_n^{k_n})D(z_1^{l_1}\cdots z_n^{l_n})-\s(z_1^{l_1}\cdots
z_n^{l_n})D(z_1^{k_1}\cdots
    z_n^{k_n})=\\
    -Qq_{z_1}^{k_1}\cdots q_{z_n}^{k_n}z_1^{\alpha_1(\kbf)}\cdots
z_n^{\alpha_n(\kbf)}\\
   \cdot z_1^{l_1-G_1}\cdots z_n^{l_n-G_n}\Big (q_{z_1}^{l_1}\cdots
q_{z_n}^{l_n}z_1^{\alpha_1(\lbf)-l_1}\cdots
       z_n^{\alpha_n(\lbf)-l_n}-1\Big )+\\
    +Qq_{z_1}^{l_1}\cdots q_{z_n}^{l_n}z_1^{\alpha_1(\lbf)}\cdots
z_n^{\alpha_n(\lbf)}\cdot\\
    \cdot z_1^{k_1-G_1}\cdots z_n^{k_n-G_n}\Big (q_{z_1}^{k_1}\cdots
q_{z_n}^{k_n}z_1^{\alpha_1(\kbf)-k_1}\cdots
        z_n^{\alpha_n(\kbf)-k_n}-1\Big )=
\end{multline*}
\begin{multline*}
    =-Qq_{z_1}^{k_1+l_1}\cdots
q_{z_n}^{k_n+l_n}z_1^{\alpha_1(\kbf)+\alpha_1(\lbf)-G_1}\cdots
    z_n^{\alpha_n(\kbf)+\alpha_n(\lbf)-G_n}+\\
    +Qq_{z_1}^{k_1}\cdots q_{z_n}^{k_n}z_1^{\alpha_1(\kbf)+l_1-G_1}\cdots
z_n^{\alpha_n(\kbf)+l_n-G_n}+\\
    +Qq_{z_1}^{k_1+l_1}\cdots
q_{z_n}^{k_n+l_n}z_1^{\alpha_1(\kbf)+\alpha_1(\lbf)-G_1}\cdots
    z_n^{\alpha_n(\kbf)+\alpha_n(\lbf)-G_n}-\\
    -Qq_{z_1}^{l_1}\cdots q_{z_n}^{l_n}z_1^{\alpha_1(\lbf)+k_1-G_1}\cdots
z_n^{\alpha_n(\lbf)+k_n-G_n}=\\
    =Qq_{z_1}^{k_1}\cdots q_{z_n}^{k_n}z_1^{\alpha_1(\kbf)+l_1-G_1}\cdots
z_n^{\alpha_n(\kbf)+l_n-G_n}-\\
    -Qq_{z_1}^{l_1}\cdots q_{z_n}^{l_n}z_1^{\alpha_1(\lbf)+k_1-G_1}\cdots
z_n^{\alpha_n(\lbf)+k_n-G_n}.
\end{multline*}This gives us
\begin{align*}
    [d_{\kbf},d_{\lbf}]_\s&=
        Qq_{z_1}^{l_1}\cdots q_{z_n}^{l_n}d_{\alpha_1(\lbf)+k_1-G_1,\dots,
\alpha_n(\lbf)+k_n-G_n}-\\
        &\quad-Qq_{z_1}^{k_1}\cdots
q_{z_n}^{k_n}d_{\alpha_1(\kbf)+l_1-G_1,\dots,\alpha_n(\kbf)+l_n-G_n}.
\end{align*}
Using the $\d$ we calculated before we can deduce a deformed
Jabobi identity as
\begin{align*}
    \cs_{k,l,h}\Big (
    &q_{z_1}^{k_1}\cdots q_{z_n}^{k_n}\big
[d_{\alpha_1(\kbf),\dots,\alpha_n(\kbf)},[d_{\lbf},d_{\hbf}]_\s
\big
    ]_\s+\\
    &+q_{z_1}^{G_1}\cdots q_{z_n}^{G_n}z_1^{\d_1}\cdots z_n^{\d_n} \big
[d_{\kbf},[d_{\lbf},d_{\hbf}]_\s \big ]_\s
    \Big )=0
\end{align*}from which we see that we get a hom-Lie algebra if all $\d_k$ are zero,
that is, if
\begin{align*}
    \begin{bmatrix}
        S_{1,1}-1 & S_{1,2} & \dots& S_{1,n}\\
        S_{2,1}& S_{2,2}-1 &\dots & S_{2,n}\\
        \vdots  &\dots\\
        S_{n,1}& & \dots &S_{n,n}-1
    \end{bmatrix}^T
    \begin{bmatrix}
        G_1\\
        G_2\\
        \vdots\\
        G_n
    \end{bmatrix}=0
\end{align*} which means that
\begin{align*}
    \ker\begin{bmatrix}
        S_{1,1}-1 & S_{1,2} & \dots& S_{1,n}\\
        S_{2,1}& S_{2,2}-1 &\dots & S_{2,n}\\
        \vdots  &\dots\\
        S_{n,1}& & \dots &S_{n,n}-1
    \end{bmatrix}^T \bigcap\, \Z^n\neq \emptyset.
\end{align*} In this case we get the Jacobi identity
\begin{align*}
        \cs_{k,l,h}\Big (
    &q_{z_1}^{k_1}\cdots q_{z_n}^{k_n}\big
    [d_{\alpha_1(\kbf),\dots,\alpha_n(\kbf)},[d_{\lbf},d_{\hbf}]_\s
    \big ]_\s+ q_{z_1}^{G_1}\cdots q_{z_n}^{G_n}\big
    [d_{\kbf},[d_{\lbf},d_{\hbf}]_\s \big ]_\s \Big )=0.
\end{align*}

We summarize the obtained results in the following theorem.
\begin{thm} \label{th:qWittmultdimmonom}
    The $\C$-linear space $$\M=
    \bigoplus_{\lbf\in\Z^n} \C\cdot d_{\lbf}$$
    spanned by $d_{\lbf}=-z_1^{l_1}\cdots z_n^{l_n}D$, where $D$ is given by
    $$D=Q\frac{\id-\s}{z_1^{G_1}\cdots z_n^{G_n}},$$ can be endowed with a bracket
    defined on generators (by (\ref{eq:GenWittProdDef})) as
    \begin{align*}
        [d_{\kbf},d_{\lbf}]_\s=q_{z_1}^{k_1}\cdots q_{z_n}^{k_n}d_{\a_1(\kbf),\dots,
\a_n(\kbf)}d_{\lbf}-
        q_{z_1}^{l_1}\cdots q_{z_n}^{l_n}d_{\a_1(\lbf),\dots, \a_n(\lbf)}d_{\kbf}
    \end{align*}
    and satisfying relations
    \begin{align*}
    [d_{\kbf},d_{\lbf}]_\s&=
        Qq_{z_1}^{l_1}\cdots q_{z_n}^{l_n}d_{\alpha_1(\lbf)+k_1-G_1,\dots,
\alpha_n(\lbf)+k_n-G_n}-\\
        &\quad-Qq_{z_1}^{k_1}\cdots
        q_{z_n}^{k_n}d_{\alpha_1(\kbf)+l_1-G_1,\dots,\alpha_n(\kbf)+l_n-G_n}.
    \end{align*} The bracket satisfies the $\s$-deformed Jacobi-identity
    \begin{align*}
        \cs_{k,l,h}\Big (
        &q_{z_1}^{k_1}\cdots q_{z_n}^{k_n}\big
[d_{\alpha_1(\kbf),\dots,\alpha_n(\kbf)},[d_{\lbf},d_{\hbf}]_\s
\big
        ]_\s+\\
        &+q_{z_1}^{G_1}\cdots q_{z_n}^{G_n}z_1^{\d_1}\cdots z_n^{\d_n} \big
[d_{\kbf},[d_{\lbf},d_{\hbf}]_\s \big ]_\s
        \Big )=0.
    \end{align*}
    Furthermore, a hom-Lie algebra is obtained if the eigenvalue problem
    $$S\cdot G=G$$ (where ''$\cdot$'' means product of matrices) has a solution
$G\in\Z^n$. We then get the Jacobi identity
    \begin{align*}
        \cs_{k,l,h}\Big (
        &q_{z_1}^{k_1}\cdots q_{z_n}^{k_n}\big
[d_{\alpha_1(\kbf),\dots,\alpha_n(\kbf)},[d_{\lbf},d_{\hbf}]_\s
\big
        ]_\s+ q_{z_1}^{G_1}\cdots q_{z_n}^{G_n}\big [d_{\kbf},[d_{\lbf},d_{\hbf}]_\s
\big ]_\s \Big )=0.
    \end{align*}
\end{thm}

\section{A deformation of the Virasoro algebra} \label{sec:cextqVirqWittlinear}
In this subsection we will prove existence and uniqueness (up to
isomorphism of hom-Lie algebras) of a one-dimensional central
extension of the hom-Lie algebra $(\mfD_\s(\C[t,t^{-1}]),\vs)$
constructed in section Section \ref{sec:qwitt}, in the case when
$q$ is not a root of unity. The obtained hom-Lie algebra is a
$q$-deformation of the Virasoro algebra.

\subsection{Uniqueness of the extension}
Let $\A=\C [t,t^{-1}]$, $\s$ be the algebra endomorphism on $\A$
satisfying $\s(t)=qt$, where $0,1\neq q\in\C$ is not a root of
unity, and set $L=\mfD_{\s}(\A)$. Then $L$ can be given the
structure of a hom-Lie algebra $(L,\vs)$ as described in Section
\ref{sec:qwitt}.

Let
$$\begin{CD} 0@>>>(\C,\id_\C)@>\iota>> (\hat L,\hat\vs)
@>\pr>>(L,\vs)@>>> 0
\end{CD}$$
be a short exact sequence of hom-Lie algebras and hom-Lie algebra
homomorphisms. In other words, let $(\hat L,\hat\vs)$ be a
one-dimensional central extension of $(L,\vs)$ by $(\C,\id_\C)$.
We also set $\c=\iota(1)$.

Choose a linear section $s:L\to\hat L$ and let
$g_s\in\mathrm{Alt}^2(L,\C)$ be the corresponding ''2-cocycle'' so
that (\ref{lift_brack}) is satisfied for $x,y\in L$. Let $\{d_n\}$
denote the basis (\ref{eq:qwittbasis}) of $L$. Define a linear map
$s':L\to\hat L$ by
\begin{align*}
    s'(d_n)&=\left\{\begin{array}{ll}
    s(d_n) &\text{if } n=0\\
    s(d_n)-\frac{1}{\{n\}}\iota\circ g_s(d_0,d_n)\c &\text{if } n\neq 0
\end{array}\right.
\end{align*}
Then $s'$ is also a section. Using the calculation
(\ref{eq:stildeg}) and the commutation relation
(\ref{eq:qwittcomrel}) we get
\begin{align*}
    \iota\circ g_{s'}(d_m,d_n)&=\iota\circ g_s(d_m,d_n)+(s-s')([d_m,d_n]_L)=\\
    &=\left\{\begin{array}{ll}
    \iota\circ g_s(d_m,d_n) &\text{if } m+n=0\\
    \iota\circ g_s(d_m,d_n)+\frac{\{m\}-\{n\}}{\{m+n\}}\iota\circ
    g_s(d_0,d_{m+n}) &\text{if } m+n\neq 0
    \end{array}\right. .
\end{align*}
In particular we have $g_{s'}(d_0,d_n)=0$ for any $n\in\Z$.
According to the calculations in Section \ref{sec:extofhlas}, the
''2-cocycle'' $g_{s'}$ must satisfy (\ref{eq:twococycle}) for any
$a,b,c\in L$. Thus we have,
$$\cs_{k,l,m}g_{s'}\big((\id+\vs)(d_k),[d_l,d_m]_L\big)=0$$
for $k,l,m\in\Z$. Substituting the definition
(\ref{eq:hOnBasisElement}) of $\vs$ and using the commutation
relation (\ref{eq:qwittcomrel}) again we get
\begin{align}\label{eq:acond}
    \cs_{k,l,m}(1+q^k)\big(\{l\}-\{m\}\big)a(k,l+m)=0,
\end{align}
 where we for simplicity have put $a(m,n)=g_{s'}(d_m,d_n)$ for $m,n\in\Z$.
 Using (\ref{eq:acond}) with $k=0$ and that $a(0,n)=0$ for any $n\in\Z$ we obtain
$$(1+q^l)\{m\}a(l,m)+(1+q^m)(-\{l\})a(m,l)=0,$$
or, since $a$ is alternating,
$$\big((1+q^l)\frac{q^m-1}{q-1}+(1+q^m)\frac{q^l-1}{q-1})\big)a(l,m)=0,$$
which simplifies to
$$2\frac{q^{l+m}-1}{q-1}a(l,m)=0.$$
This shows that $a(l,m)=0$ unless $l+m=0$. Setting $b(m)=a(m,-m)$
we have so far
$$[s'(d_m),s'(d_n)]_{\hat L}=(\{m\}-\{n\})s'(d_{m+n})+\d_{m+n,0}b(m)\c.$$
Use (\ref{eq:acond}) with $k=-n-1$, $l=n$, $m=1$ we get
\begin{align*}
    (1+q^{-n-1})\frac{q^n-q}{q-1}a(-n-1,n+1)&+
    (1+q^n)\frac{q-q^{-n-1}}{q-1}a(n,-n)+\\
    &+(1+q)\frac{q^{-n-1}-q^n}{q-1}a(1,-1)=0,
\end{align*}
or, after multiplication by $q^{n+1}$,
\begin{align} \label{eq:recurrence1}
    q(1+q^{n+1})\{n-1\}b(n+1)=(1+q^n)\{n+2\}b(n)-(1+q)\{2n+1\}b(1).
\end{align} This is a second order linear recurrence equation in $b$.
\begin{lem} \label{lem:nonrootof1} The functions $b_1,b_2:\Z\to\C$ defined by
\begin{align*}
    b_1(m)&=\frac{q^{-m}}{1+q^{m}}\{m-1\}\{m\}\{m+1\}\\
    b_2(m)&=q^{-m}\{2m\}
\end{align*}
are two linear independent solutions of (\ref{eq:recurrence1}).
\end{lem}
\begin{proof}
Substituting $b_1$ for $b$ in (\ref{eq:recurrence1}) the left hand
side equals
$$q\{n-1\}q^{-n-1}\{n\}\{n+1\}\{n+2\},$$
while the right hand side becomes
$$\{n+2\}q^{-n}\{n-1\}\{n\}\{n+1\}-0.$$
These expressions are equal. To prove that $b_2$ is also a
solution requires some calculations:
\begin{multline*}
    q\{n-1\}(1+q^{-1})\{2n+2\}-\{n+2\}(1+q^{-n})\{2n\}+\{2n+1\}(1+q^{-1})\{2\}=\\
    =(\{n\}-q^{n-1})(q+q^{-n})(\{2n\}+q^{2n}+q^{2n+1})-\\
    -(\{n\}+q^n(q+1))(1+q^{-n})\{2n\}+(\{2n\}+q^{2n})(1+q^{-1})(1+q)=\\
    =\{2n\}\Big(\{n\}(q+q^{-n})-q^n-q^{-1}-\{n\}(1+q^{-n})-\\
    -(q+1)q^n-1-q+2+q+q^{-1}\Big)+\\
    +(q+1)q^{2n}\Big(q\{n\}+q^{-n}\{n\}-q^n-q^{-1}+1+q^{-1}\Big)=\\
    =\{2n\}\Big(q^n-1-q^n-q^n(q+1)-1+2\Big)+q^{2n}\frac{q+1}{q-1}(-q^n+q^n)=\\
    =\frac{q^{2n}-1}{q-1}(q+1)(-q^n)+q^{2n}\frac{q+1}{q-1}(-q^n+q^n)=\\
    =\frac{q+1}{q-1}(q^n-q^{2n-n})=0.
\end{multline*}
It remains to show that $b_1$ and $b_2$ are linear independent. If
$$\lambda b_1+\mu b_2=0,$$
then evaluation at $m=1$ gives $\mu q^{-1}(1+q)=0$ so $\mu=0$.
Since $b_1$ is nonzero, we must have $\lambda=0$ also.
\end{proof}
Thus we have
$$b(m)=\alpha b_1(m)+\beta b_2(m)$$
for some $\alpha,\beta\in\C$. In terms of $g_{s'}$ this means that
$$g_{s'}(d_m,d_n)=\d_{m+n,0}(\alpha b_1(m)+\beta b_2(m)).$$
Define now yet another section $s'':L\to\hat L$ by
$$s''(d_m)=s'(d_m)+\delta_{m,0}\beta \c.$$
Then
\begin{align*}
    \iota\circ g_{s''}(d_m,d_n)&=\iota\circ g_{s'}(d_m,d_n)+(s'-s'')([d_m,d_n]_L)=\\
    &=\d_{m+n,0}(\alpha b_1(m)+\beta
    b_2(m))\c-(\{m\}-\{n\})\d_{m+n,0}\beta\c. \\
    &=\d_{m+n,0}(\alpha b_1(m)+\beta b_2(m)-\beta q^{-m}\{2m\})\c=\\
    &=\d_{m+n,0}\alpha b_1(m)\c,
\end{align*}
where we used that $\{m\}-\{-m\}=q^{-m}\{2m\}$. If $\alpha=0$ we
have a trivial extension. Otherwise we set $\c'=6\alpha\c$.

It remains to determine the homomorphism $\hat\vs$. Using equation
(\ref{eq:sigmahattfs}) we have for $x\in\hat L$,
$$\hat\vs(x)=s''\circ\vs\circ\pr(x)+\iota\circ f_{s''}(x)$$
for some linear function $f_{s''}:\hat L\to\C$. To determine
$f_{s''}$, first use (\ref{eq:sigmaafrakfs}):
$$f_{s''}(\c')=f_{s''}(\iota(6\a))=\id_\C(6\a)=6\a.$$
Hence
$$\hat\vs(\c')=\c'.$$
Next, we use (\ref{eq:bothgsandfs}) in Theorem \ref{thm:homext1}
to get
$$f_{s''}\big([s''(d_m),s''(d_n)]_{\hat L}\big)=g_{s''}(\vs(d_m),\vs(d_n)).$$ By
(\ref{lift_brack}) we see that
\begin{align*}
    [s''(d_m),s''(d_n)]_{\hat{L}} &=s''([d_m,d_n]_L)+\iota\circ g_{s''}(d_m,d_n)=\\
    &=\big (\{m\}-\{n\}\big )s''(d_{m+n})+\iota\circ g_{s''}(d_m,d_n)
\end{align*} and so
$$\big(\{m\}-\{n\}\big)f_{s''}\big(s''(d_{m+n})\big)+f_{s''}\big(\iota\circ
g_{s''}(d_m,d_n))=q^{m+n}g_{s''}(d_m,d_n)$$ which is equivalent to
$$\big(\{m\}-\{n\}\big)f_{s''}\big(s''(d_{m+n})\big)=(q^{m+n}-1)g_{s''}(d_m,d_n)$$
for all integers $m,n$. But the right hand side is identically
zero ($g_{s''}$ being a multiple of $\d_{m+n,0}$). Hence, taking
$m\neq 0$ and $n=0$ we get $f_{s''}(s''(d_m))=0$ for all nonzero
$m$. But if we take $m=1$ and $n=-1$ we also get
$f_{s''}(s''(d_0))=0$ because $\{1\}-\{-1\}=1+q^{-1}\neq 0$ since
$q$ is not a root of unity.

Putting $\hat{L}\ni L_n:=s''(d_n)$ we have proved the following
theorem.
\begin{thm}\label{thm:uniq} Every nontrivial one-dimensional central extension
of the hom-Lie algebra $\big (\mfD_\s(\A),\vs\big )$, where
$\A=\C[t,t^{-1}]$, is isomorphic to the hom-Lie algebra
$\Vir_q=(\hat{L},\hat\vs)$, where $\hat{L}$ is the non-associative
algebra with basis $\{L_n:n\in\Z\}\cup\{\c\}$ and relations
\begin{align*}
    [\c,\hat{L}]_{\hat{L}}&=0,\\
    [L_m,L_n]_{\hat{L}}&=\big(\{m\}-\{n\}\big)L_{m+n}+\d_{m+n,0}\frac{q^{-m}}{6(1+q^m)}\{m-1\}\{m\}\{m+1\}\c,
\end{align*}
and $\hat\vs:\hat{L}\to\hat{L}$ is the endomorphism of $\hat{L}$
defined by
\begin{align*}
    \hat\vs(L_n)&=q^nL_n, \\
    \hat\vs(\c)&=\c.
\end{align*}
\end{thm}

\subsection{Existence of a non-trivial extension}
We now proceed to prove the following result. Let
$\A=\C[t,t^{-1}]$.
\begin{thm} There exists a non-trivial central extension of
$\big (\mfD_\s(\A),\vs\big )$ by $(\C,\id_\C)$.
\end{thm}
\begin{proof}
We set $L:=\mfD_\s(\A)$ for brevity and define $g:L\times
L\longrightarrow\C$ by setting
$$g(d_m,d_n):=\d_{m+n,0}\frac{q^{-m}}{1+q^m}\{m-1\}\{m\}\{m+1\}, \quad\text{for }
m,n\in\Z$$ and extending using the bilinearity. We also define a
linear map $f:L\oplus\C\longrightarrow\C$ by
$$f(x,a)=a\quad\text{for } x\in L,\, a\in\C.$$ Our goal is to use Theorem
\ref{thm:homext2} which means that we have to verify that $g$ and
$f$ satisfy the necessary conditions. First of all, using that
$\{-n\}=-q^{-n}\{n\}$, we note
\begin{multline*}
    g(d_m,d_n)=\d_{m+n,0}\frac{q^{-m}}{1+q^m}\{m-1\}\{m\}\{m+1\}=\\
    =\d_{n+m,0}\frac{q^n}{1+q^{-n}}\{-n-1\}\{-n\}\{-n+1\}=\\
    =-\d_{n+m,0}\frac{q^{-n}}{1+q^n}\{n-1\}\{n\}\{n+1\}=-g(d_n,d_m).
\end{multline*}This shows that $g$ is alternating. That (\ref{eq:cfvsa}) holds is
immediate. To check (\ref{eq:cbothgsandfs}) let $m,n\in\Z$. Then
\begin{multline*}
    g(\vs(d_m),\vs(d_n))=g(q^md_m,q^nd_n)=q^{m+n}g(d_m,d_n)=\\
    =q^{m+n}\d_{m+n,0}\frac{q^{-m}}{1+q^m}\{m-1\}\{m\}\{m+1\}=\\
    =g(d_m,d_n)=f([d_m,d_n]_L,g(d_m,d_n))
\end{multline*} It remains to verify (\ref{eq:ctwococycle2}). By trilinearity it is
enough to assume that $(x,y,z)=(d_k,d_l,d_m)$ for some
$k,l,m\in\Z$. Moreover, if $k+m+l\neq 0$ then
(\ref{eq:ctwococycle2}) holds trivially due to the Kronecker delta
in the definition of $g$. Thus we can assume $k+m+l=0$. We then
have
\begin{multline*}
    \cs_{k,l,m} g\big ((\id+\vs)d_k,[d_l,d_m]_L\big )=\cs_{k,l,m} (1+q^k)\big
(\{l\}-\{m\}\big
    )g(d_k,d_{l+m})=\\
    =\big(\{l\}-\{-k-l\}\big)q^{-k}\{k-1\}\{k\}\{k+1\}+\\
    \quad+\big(\{-k-l\}-\{k\}\big)q^{-l}\{l-1\}\{l\}\{l+1\}+\\
    \quad+\big(\{k\}-\{l\}\big)q^{k+l}\{-k-l-1\}\{-k-l\}\{-k-l+1\}=\\
    =-\{-k-l\}\Big(q^{-k}\{k-1\}\{k\}\{k+1\}-q^{-l}\{l-1\}\{l\}\{l+1\}-\\
    \quad-(\{k\}-\{l\})q^{k+l}\{-k-l-1\}\{-k-l+1\}\Big)+\\
    \quad+q^{-k-l}\{k\}\{l\}\Big(q^l\{k-1\}\{k+1\}-q^k\{l-1\}\{l+1\}\Big).
\end{multline*} The second factor in the first term equals
\begin{multline*}
    q^{-k}\{k-1\}\{k\}\{k+1\}-q^{-l}\{l-1\}\{l\}\{l+1\}-\\
    -(\{k\}-\{l\})q^{k+l}\{-k-l-1\}\{-k-l+1\}=\\
    =q^{-k}(\{k\}-q^{k-1})\{k\}(\{k\}+q^k)-q^{-l}(\{k\}-q^{l-1})\{l\}(\{l\}+q^l)-\\
    -(\{k\}-\{l\})q^{k+l}q^{-2k-2l}(\{k+l\}+q^{k+l})(\{k+l\}-q^{k+l-1})=\\
    =q^{-k}\{k\}^3+(1-q^{-1})\{k\}^2-q^{k-1}\{k\}-q^{-l}\{l\}^3-(1-q^{-1})\{l\}^2+q^{l-1}\{l\}-\\
    -(\{k\}-\{l\})\Big(q^{-k-l}(\{k\}+q^k\{l\})(\{l\}+q^l\{k\})+(1-q^{-1})\{k+l\}-q^{k+l-1}\Big)=\\
    =(1-q^{-1})\{k\}^2-q^{k-1}\{k\}-(1-q^{-1})\{l\}^2+q^{l-1}\{l\}-\\
    -q^{-k-l}\{k\}\{l\}\Big((1+q^{k+l})(\{k\}-\{l\})+q^k\{l\}-q^l\{k\}\Big)-\\
    -(\{k\}-\{l\})\Big((1-q^{-1})\{k+l\}-q^{k+l-1}\Big)=\\
    =\{k\}\Big((1-q^{-1})(\{k\}-\{k+l\})-q^{k-1}+q^{k+l-1}\Big)-\\
    -\{l\}\Big((1-q^{-1})(\{l\}-\{k+l\})-q^{l-1}+q^{k+l-1}\Big)-\\
    -q^{-k-l}\{k\}\{l\}\Big((1+q^{k+l})(\{k\}-\{l\})+q^k\{l\}-q^l\{k\}\Big).
\end{multline*}
The first of the three terms in the last equality above is equal
to
\begin{multline*}
    (1-q^{-1})(\{k\}-\{k+l\})-q^{k-1}+q^{k+l-1}=\\
    =q^{-1}(q-1)\frac{q^k-1-q^{k+l}+1}{q-1}-q^{k-1}+q^{k+l-1}=0.
\end{multline*}
Similarly the second term vanishes. Using that
$\{a+b\}=\{a\}+q^a\{b\}$ we see that the whole expression is equal
to
\begin{multline*}
    q^{-k-l}\{k\}\{l\}\Big(q^l\{k-1\}\{k+1\}-q^k\{l-1\}\{l+1\}+\\
    +\{-k-l\}(
    (1+q^{k+l})(\{k\}-\{l\})+q^k\{l\}-q^l\{k\})\Big)=\\
    =q^{-k-l}\{k\}\{l\}\Big(q^l(\{k\}^2+q^k(1-q^{-1})\{k\}-q^{2k-1})-\\
    -q^k(\{l\}^2+q^l(1-q^{-1})\{l\}-q^{2l-1})
    -\{k+l\}(\{k\}-\{l\})\Big)=\\
    =q^{-k-l}\{k\}\{l\}\Big(q^l\{k\}^2-q^{k+l+1}-q^k\{l\}^2+q^{k+l+1}
    -\{k+l\}(\{k\}-\{l\})\Big)=\\
    =q^{-k-l}\{k\}\{l\}\{k+l\}\Big(q^l\{k\}-q^k\{l\}-\{k\}+\{l\}\Big)=\\
    =\{k\}\{l\}\{k+l\}\Big (\{l+k\}-\{k+l\}\Big )=0.
\end{multline*}Thus (\ref{eq:ctwococycle2}) holds for all $x,y,z\in L$. Hence, by
Theorem \ref{thm:homext2}, there exists a central extension
of $(L,\vs)$ by $(\afrak,\vs_{\afrak})$. \\

Suppose this extension is trivial. Then by Remark
\ref{rem:trivial} there is a linear map $s_1$ such that
$$g(d_m,d_n)=s_1([d_m,d_n]_L),$$ or
$$\d_{m+n,0}\frac{q^{-m}}{1+q^m}\{m-1\}\{m\}\{m+1\}=(\{m\}-\{n\})s_1(d_{m+n}),$$ for
$m,n\in\Z$. Taking $m=1$ and $n=-1$ gives $s_1(d_0)=0$. On the
other hand, setting $m=2$ and $n=-2$ yields $s_1(d_0)\neq 0$. This
contradiction shows that the extension is non-trivial.
\end{proof}
\begin{remark}
    The coefficient in the central extension part in Theorem \ref{thm:uniq} is
$1/6\cdot g(d_m,d_n)$, where
    $g(d_m,d_n)$
     is from the above Theorem. This factor $1/6$ is easily obtained by rescaling
$\c$. The reason for this factor
     in Theorem \ref{thm:uniq} is that for the classical undeformed Virasoro algebra
one usually
     rescales by a factor $1/12$ in the central extension term. Now by taking $q=1$
in Theorem \ref{thm:uniq},
     we thus get the classical undeformed Virasoro algebra including the usually
chosen scaling factor $1/12$.
\end{remark}

\begin{remark}
If $\t$ is an automorphism of $\A$ and $\Delta$ is a
$(\s,\t)$-derivation on $A$ we can still define a product on
$\A\cdot\Delta$ by
$$[a\cdot\Delta,b\cdot\Delta]_{\s,\t,\Delta}= (\s(a) \cdot \Delta) \circ (b \cdot \Delta) -
(\s(b) \cdot \Delta) \circ (a \cdot \Delta) =
\big(\s(a)\Delta(b)-\s(b)\Delta(a)\big)\cdot\Delta.$$ For example,
if we take $\A=\C[t,t^{-1}]$, $\s(t)=qt$ and $\t=\s^{-1}$, and the
symmetric $q$-difference operator $$\Delta=\frac{\t-\s}{q^{-1}-q}
: f(t) \mapsto \frac{f(q^{-1}t)-f(qt)}{q^{-1}-q},$$ then our
bracket will be
\begin{equation}
[d_n,d_m]= q^n d_n d_m - q^m d_m d_n =
\frac{q^{n-m}-q^{m-n}}{q-q^{-1}}d_{n+m}=[n-m]d_{n+m},\label{eq:satoshogerled}
\end{equation}
where $d_n=-t^n\cdot\Delta$ and $[k]=(q^k-q^{-k})/(q-q^{-1})$ is
the symmetric $q$-number. This is easily calculated by direct
substitution. The right hand side of this commutation relation
\eqref{eq:satoshogerled} coincides with the right hand side of the
defining relations (1) in the $q$-deformation of Witt algebra
considered in \cite{AizSato1,CurtrZachos1}, but the left hand
bracket side turns out to be slightly different. When
$q\rightarrow 1$ the defining relations for the classical Witt
algebra are recovered in both cases.
\end{remark}


\begin{thebibliography}{99}
\bibitem{AizSato1} Aizawa, N., Sato, H.-T., \emph{$q$-deformation of the
Virasoro algebra with central extension}, Phys. Lett. B, 256
(1999), no. 2, 185--190.
\bibitem{AvanFrappatRossiSorba} Avan, J., Frappat, L., Rossi,
M., Sorba, P., \emph{Central extensions of classical and quantum
$q$-Virasoro algebras}, Phys. Lett. A 251 (1999), no. 1, 13--24.
%
\bibitem{BelovChalt1} Belov, A. A., Chaltikian, K. D., \emph{$q$-deformation
of Virasoro algebra and lattice conformal theories}, Modern Phys.
Lett. A 8 (1993), no. 13, 1233--1242.
\bibitem{BelovChalt2} Belov, A. A., Chaltikian, K. D., \emph{$q$-deformation of
Virasoro algebra and lattice conformal theories}, Zap. Nauchn.
Sem. S.-Peterburg. Otdel. Mat. Inst. Steklov. (POMI) 235 (1996),
Differ. Geom. Gruppy Li i Mekh. 15-2, 217--227, 306; translation
in J. Math. Sci. (New York) 94 (1999), no. 4, 1581--1588
\bibitem{BlochS1} Bloch S., \emph{Zeta values and differential operators on the
circle}, J. Algebra, 182 (1996), 476--500.
\bibitem{ChaiIsLukPopPres} Chaichian, M., Isaev, A. P.,
Lukierski, J., Popowicz, Z., Pre\v snajder,~P.,
\emph{$q$-deformations of Virasoro algebra and conformal
dimensions}, Phys. Lett. B 262 (1991), no. 1, 32-38.
\bibitem{ChaiKuLuk} Chaichian, M., Kulish, P.,
Lukierski, J., \emph{$q$-deformed Jacobi identity, $q$-oscillators
and $q$-deformed infinite-dimensional algebras}, Phys. Lett. B 237
(1990), no. 3-4, 401--406.
\bibitem{ChaPopPres} Chaichian, M., Popowicz, Z., Pre\v
snajder, P., \emph{$q$-Virasoro algebra and its relation to the
$q$-deformed KdV system}, Phys. Lett. B 249 (1990), no. 1, 63--65.
\bibitem{ChaiPres1} Chaichian, M., Pre\v snajder, P., \emph{Sugawara construction
and the $q$-deformation of Virasoro algebra}, in Quantum groups
and related topics (Wroc{\l}aw, 1991), 3--12, Math. Phys. Stud.,
13, Kluwer Acad. Publ., Dordrecht, 1992.
\bibitem{ChaiPres2} Chaichian, M., Pre\v snajder, P., \emph{On the $q$-Sugawara
construction for the Virasoro (super) algebra}, in Quantum
symmetries (Clausthal, 1991), 352--365, World Sci. Publishing,
River Edge, NJ, 1993.
\bibitem{ChaiPres3} Chaichian, M., Pre\v snajder, P.,
\emph{$q$-Virasoro algebra, $q$-conformal dimensions and free
$q$-superstring}, Nuclear Phys. B 482 (1996), no. 1-2, 466--478.
\bibitem{ChakJag} Chakrabarti, R., Jagannathan, R., \emph{A
$(p,q)$-deformed Virasoro algebra}, J. Phys. A 25 (1992), no. 9,
2607--2614.
\bibitem{ChungWS} Chung, W.-S., \emph{Two parameter deformation of Virasoro
algebra}, J. Math. Phys. 35 (1994), no. 5, 2490--2496.
\bibitem{CurtrZachos1} Curtright, T. L., Zachos, C. K.,
\emph{Deforming maps for quantum algebras}, Phys. Lett. B 243
(1990), no. 3, 237--244.
\bibitem{DevSav} Devchand, Ch., Saveliev, M. V.,
\emph{Comultiplication for quantum deformations of the centreless
Virasoro algebra in the continuum formulation}, Phys. Lett. B 258
(1991), no. 3-4, 364--368.
\bibitem{FairNuyZach} Fairlie, D. B., Nuyts, J., Zachos, C. K., \emph{A presentation
for the Virasoro and Super-Virasoro algebras}, Commun. Math. Phys.
117 (1988), 595-614.
\bibitem{HelSil-book}
{\rm L. Hellstr{\"o}m, S. D. Silvestrov,} \emph{Commuting Elements
in $q$-Deformed Heisenberg Algebras},
World Scientific, 2000, 256 pp.\\
 (ISBN: 981-02-4403-7).
\bibitem{HuN} Hu, N., \emph{Quantum group structure of the
$q$-deformed Virasoro algebra}, Lett. Math. Phys. 44 (1998), no.
2, 99--103.
\bibitem{JellalSatoHT} Jellal A., Sato, H.-T., \emph{FFZ realization of the
deformed super Virasoro algebra -- Chaichian-Pre\v{s}naider type},
Phys. Lett. B 483 (2000),  451--455.
\bibitem{KacRadul} Kac, V., Radul, A., \emph{Quasifinite highest weight modules over
the Lie algebra of differential operators on the circle}, Commun.
Math. Phys. 157 (1993), 429-457.
\bibitem{Kassel1} Kassel, C., \emph{Cyclic homology of differential operators,
the Virasoro algebra and a $q$-analogue}, Commun. Math. Phys. 146
(1992), 343-351.
\bibitem{KemmokuSatoHT} Kemmoku, R., Sato, H.-T., \emph{Deformed fields and Moyal
construction of deformed super Virasoro algebra}, Nucl. Phys. B
595 (2001), 689-709.
\bibitem{KirkmanProcesiSmall} Kirkman, E., Procesi, C., Small, L.,
\emph{A $q$-analogue for Virasoro algebra}, Comm. Alg. 22 (1994),
3755-3774.
\bibitem{KhesinLyubRoger} Khesin, B., Lyubashenko, V., Roger, C.,
\emph{Extensions and contractions of Lie algebra of
$q$-Pseudodifferential symbols on the circle}, J. Func. Anal. 143
(1997), 55-97.
\bibitem{Li} Li, W.-L., \emph{$2$-Cocycles on the algebra of differential operators},
J. Alg. 122 (1989), 64-80.
\bibitem{LiuKQ1} Liu, K. Q., \emph{Some results on $q$-deformations of the Virasoro
algebra}, in Proc. Conf. on Quantum Topology (Manhattan, KS,
1993), World Sci. Publ., River Edge, NJ, 1994, 259--268.
\bibitem{LiuKQ2} Liu, K. Q., \emph{Indecomposable representations of the
$q$-deformed Virasoro algebra}, Math. Z. 217 (1994), no. 1,
15--35.
\bibitem{LiuKQ3} Liu, K. Q., \emph{A class of Harish-Chandra modules for the
$q$-deformed Virasoro algebra}, J. Algebra, 171 (1995), no. 2,
606--630.
\bibitem{Mazorchuck1} Mazorchuck, V.,
\emph{On simple modules over $q$-analog for the Virasoro algebra},
Hadronic J. 21 (1998), no. 5, 541--550.
\bibitem{MebAisBouMaa} Mebarki, N., Aissaoui, H., Boudine,
A., Maasmi, A., \emph{The $q$-deformed Virasoro algebra},
Czechoslovak J. Phys. 47 (1997), no. 8, 755--759.
\bibitem{Narg-Quij} Narganes-Quijano, F. J., \emph{Cyclic
representations of a $q$-deformation of the Virasoro algebra}, J.
Phys. A 24 (1991), no. 3, 593--601.
\bibitem{OsbornPassman} Osborn, J. M., Passman, D. S.,
\emph{Derivations of skew-polynomial rings}, J. Algebra, 176
(1995), 417--448.
\bibitem{OsbornZhao} Osborn, J. M., Zhao, K.,
\emph{A characterization of the Block Lie algebra and its
$q$-forms in characteristic $0$}, J. Algebra, 207 (1998), no.2,
367--408.
\bibitem{Polychronakos} Polychronakos, A. P., \emph{Consistency conditions
and representations of a $q$-deformed Virasoro algebra}, Phys.
Lett. B 256 (1991), no. 1, 35--40.
\bibitem{SatoHT1} Sato, H.-T., \emph{Realizations of $q$-deformed Virasoro algebra},
Progress of Theoretical Physics, 89 (1993), no. 2, 531-544.
\bibitem{SatoHT2} Sato, H.-T., \emph{$q$-Virasoro operators from an
analogue of the Noether currents}, Z. Phys. C 70 (1996), no. 2,
349--355.
\bibitem{SatoHT3} Sato, H.-T., \emph{OPE formulae for deformed super-Virasoro
algebras}, Nucl. Phys. B 471 (1996), 553-569.
\bibitem{SatoHT4} Sato, H.-T., \emph{Deformation of super Virasoro algebra
in non-commutative quantum superspace}, Phys. Lett. B 415 (1997),
170-174.
\bibitem{Su1} Su, Y., \emph{$2$-Cocycles on the Lie algebras of generalized
differential operators}, Comm. Alg. 30 (2) (2002), 763-782.
\bibitem{Su2} Su, Y., \emph{Classification of quasifinite modules over the Lie
algebras of Weyl type}, Adv. Math. 174 (2003), 57-68.
\bibitem{ZhaoK} Zhao, K.,
\emph{The $q$-Virasoro-like algebra}, J. Algebra, 188 (1997),
506--512.
\bibitem{ZhaCZZhaoWZ} Zha, C.-Z., Zhao, W.-Z.,
\emph{The $q$-deformation of super high-order Virasoro algebra},
J. Math. Phys., 36 (2) (1995), 967--979.
\end{thebibliography}
\end{document}